\documentclass[A4paper,11pt]{article}
\pdfoutput=1

\usepackage{latexsym}
\usepackage{epsfig}
\usepackage{verbatim}
\usepackage{multirow}
\usepackage{fancybox}
\usepackage{shadow}
\usepackage{cite}
\usepackage{amscd, amsmath, amssymb, empheq}
\usepackage{bm, dsfont}
\usepackage{graphicx}
\usepackage{newtxtext, newtxmath}
\usepackage{xcolor}
\usepackage{mdframed}
\mdfsetup{%
linecolor=white,
backgroundcolor=gray!30,}
\usepackage[caption=false]{subfig}
\allowdisplaybreaks
\usepackage{float}

\usepackage[normalem]{ulem}

\newtheorem{theorem}{\sc Theorem}[section]
\newtheorem{lemma}[theorem]{\sc Lemma}

\newcommand{\eps}{\varepsilon}

\newcommand{\proofend}{{\medskip\medskip}}
\newcommand{\proof}{{\noindent\em Proof. }}

\author{
 Bernard Chazelle
\thanks{Department of Computer Science,
       Princeton University, 
{\tt chazelle}@{\tt cs.princeton.edu }}
\and
Kritkorn Karntikoon
\thanks{Department of Computer Science,
       Princeton University, 
{\tt kritkorn}@{\tt cs.princeton.edu }}
}

\title{
The $s$-Energy and Its Applications
\thanks{This work was supported in part by NSF grant
CCF-2006125.  This article presents and unifies two nonoverlapping works
presented in preliminary form in:
B. Chazelle, K. Karntikoon,
Quick relaxation in collective motion,
Proc. 61st IEEE Conference on Decision and Control, Cancun, Mexico, 2022;
and 
B. Chazelle, K. Karntikoon,
A connectivity-sensitive approach to consensus dynamics,
SAND 2023, 2nd Symposium on Algorithmic Foundations.
}}

\date{}

\begin{document} \maketitle

\begin{abstract}

Many multi-agent systems evolve by repeatedly updating each state to a weighted average
of its neighbors, a process known as averaging dynamics, whose behavior becomes difficult 
to analyze when the interaction network varies over time. 
In recent years, the {\em $s$-energy} has emerged as a useful tool for bounding the convergence
rates of such systems, complementing the classical techniques that rely on fixed graphs. 
We derive new bounds on the $s$-energy under minimal connectivity assumptions. 
As a consequence, we obtain convergence guarantees for several models of collective dynamics and 
resolve a number of open questions in the areas. 
Our results highlight the dependence of the $s$-energy on the connectivity
of the underlying networks and use it to explain the exponential gap 
in the convergence rates of stationary and time-varying consensus systems.
\end{abstract}

\vspace{2cm}

\section{Introduction}\label{sec:intro}

We study discrete-time averaging systems defined by sequences of stochastic matrices 
whose nonzero entries specify an underlying, possibly time-varying interaction graph.
At each step, agents update their states by averaging them with those of their neighbors.
Such systems arise naturally when local interactions are encoded by networks and 
global evolution is governed by repeated application of averaging operators~\cite{
blondelHT09, bulloBk, caoMA2008, 
castellanoFL2009, charron-bost22, chatterjeeS77,
FagnaniF,  hegselmanK,
HendrickxB, HendrickxT12, jadbabaieLM03,
KLO10, Lashhab18, lorenz05, Lorenz2010,
mirtaB11,  moallemiVR06, Moreau2005,
Nedich15, nedicOOT, olshevskyT-09, tsitsiklisBA,
chazelle-Energ1-2011, chazelle-Energ2-2019,
CuckerSmale1, degroot1974, earlS, vicsekCBCS95
}.
From a mathematical perspective, these dynamics interpolate between products of 
stochastic matrices, random walks on graphs, and time-inhomogeneous Markov chains.

When the interaction graph is fixed, the evolution reduces to a random walk governed
by a single stochastic matrix. Its convergence properties are well understood through
spectral graph theory and mixing time~\cite{LPW}.
By contrast, once the graph and weights are allowed to change over time, 
classical tools based on eigenvalues or stationary distributions no longer work.
Understanding convergence in this setting therefore requires new invariants that remain meaningful 
under arbitrary graph sequences subject only to local constraints.

One such invariant is the {\em $s$-energy}, a scale-sensitive function associated with the 
evolution of averaging systems~\cite{chazelle-Energ1-2011}.
Defined independently of any particular application, the $s$-energy measures how 
pairwise distances between agents shrink over time across multiple spatial scales.
The parameter $s \in (0, 1]$ controls the scale at which distances are emphasized,
enabling the function to detect both global and local behavior.

Our main results establish new bounds on the $s$-energy under minimal assumptions on 
the interaction weights and graph structure. In particular, for systems whose interaction graph 
has at most $m$ connected components at every time step, we show that the total $s$-energy is bounded 
by a function that is polynomial in the number of agents and exponential only in $m$.
The bounds hold for arbitrary time-varying graphs, and improve further under additional assumptions.

Our result highlights the key role played by network connectivity in the $s$-energy.
It helps us resolve what has been a puzzling mystery 
for several years now: Previous work, indeed, has shown that convergence 
to equilibrium within $\eps$ is reached in time at most
$C \big(\log \frac{1}{\eps}\big)^{n-1}$, where $n$ is the number of agents
and $C$ is a parameter depending only on the networks' 
characteristics~\cite{chazelle-Energ2-2019, FagnaniF, jadbabaieLM03, olshevskyT-09, olshevTsi-13}.
When all of the communication networks are connected, however, the time bound drops to 
$C \log \frac{1}{\eps}$.  Such an exponential gap also appears in recent works on broadcast and 
consensus dynamics and hyper-torpid 
mixing~\cite{chazelleIEEE-TN20, coulouma15, El-HayekH022, El-HayekH023, WinklerPG0023}.
We explain this gap by proving that the convergence time
is actually of the form $C \big(\log \frac{1}{\eps}\big)^{m}$, where $m<n$
is the maximum number of connected components at any time.   
This result follows from our new bounds on the $s$-energy.

The abstract nature of the $s$-energy framework allows it to be applied uniformly across 
a range of models. As corollaries to our main bounds, we obtain polynomial-time convergence 
guarantees for multi-flock bird dynamics, attraction rates for opinion dynamics
with stubborn agents, an explanation of the ``Overton window'' phenomenon,
as well as bounds for several classes of randomized 
or reversible averaging processes. In each case, the system-specific analysis reduces to verifying 
that the dynamics fit within the general averaging framework.

From a broader perspective, the $s$-energy can be viewed as a nonlinear analogue of classical energies
used in the analysis of random walks and graph Laplacians. While spectral graph and Dirichlet forms
quantify contraction in fixed graphs, the $s$-energy remains well behaved under graph changes, 
adversarial updates, and non-reversible dynamics. To our knowledge, it provides the first general-purpose
tool for bounding convergence in time-varying averaging systems without imposing global regularity conditions.

The paper is organized as follows. Section~\ref{sec:prelim} introduces the averaging framework, 
formally defines the $s$-energy, and states our results mathematically.
Section~\ref{sec:proofs} proves main results: general upper and lower bounds.
Section~\ref{app1} applies the theory to several dynamical systems arising in distributed control,
opinion dynamics, and collective motion.

\section{Preliminaries and Models}\label{sec:prelim}
We now formalize the class of averaging systems studied in this paper and introduce the 
notation and assumptions used throughout. An averaging system consists of a sequence of 
stochastic matrices $(P_{t})_{t \geq 0}$, where each $P_t$ specifies a graph $G_t$ 
in the sense that $(P_t)_{ij} > 0$ if and only if $(i,j)$ is an edge of $G_t$.
The state vector $x(t) \in \mathbb{R}^n$ then evolves according to
\[
x(t+1) = P_t x(t).
\]

\noindent
Within this framework, we define the $s$-energy, a scale-sensitive functional 
associated with the evolution of the system, and recall several basic conventions that will be used throughout.

\paragraph{Averaging dynamics.}

Let $(G_t)_{t\geq 0}$ be an infinite sequence of undirected graphs
over the vertex set $[n]$. Each vertex has a self-loop.
Let $P_{t}$ be the stochastic matrix of a weighted random walk over $G_t$.
By construction, a matrix entry is positive if and only if it corresponds to an edge of $G_t$.
Each row sums up to 1 and the diagonal is positive everywhere.
We assume that the nonzero entries in $P_{t}$ are at least some fixed 
$\rho \in (0, 1/2]$, called the {\em weight threshold}.\footnote{We may 
exclude the case $\rho > 1/2$, since it corresponds
to graphs without edges (except for the self-loops).} 
Let $P_{\leq t}$ denote the product $P_t P_{t-1}\cdots P_0$.
The set of orbits $(P_{\leq t} x)_{t\geq 0}$, over all $x\in \mathbb{R}^n$,
defines an {\em averaging system}. It is sometimes called 
``consensus dynamics'' in the literature~\cite{FagnaniF}.
In addition to the general case defined above, we also consider three special cases:

\begin{itemize}
\item
{\sc rev}: \  \
a {\em reversible} averaging system assumes that
all of the individual Markov chains $P_t$ are reversible
and share the same stationary distribution.
This means that $P_t= \text{diag}(q)^{-1} M_t$, where
$M_t$ is symmetric with nonzero entries at least~1
and $q= M_t\mathbf{1}\preceq  \mathbf{1}/\rho$.
The stationary distribution is given by $q/\|q\|_1$.

\item
{\sc exp}: \  \
an {\em expanding} averaging system 
is a {\sc rev} where $q= d\mathbf{1}$ and the 
connected components of each $G_t$ are $d$-regular 
expanders. Recall that a $d$-regular expander is a graph of degree $d$
such that, for any set $X$ of at most half of the vertices,
we have $|\partial X|\geq h|X|$, where $\partial X$
is the set of edges with exactly one vertex in $X$;
the factor $h$ is called the Cheeger constant.

\item
{\sc ran}: \ 
a {\em random} averaging system is a {\sc rev} where $q= d\mathbf{1}$
with each $P_t$ is the simple random walk over
a random $(d-1)$-regular graph $G_t$.

\end{itemize}

\paragraph{The $\bm{s}$-energy.}

\begin{figure}[H]
\begin{center}
\hspace{0cm}
\includegraphics[height=3.5cm]{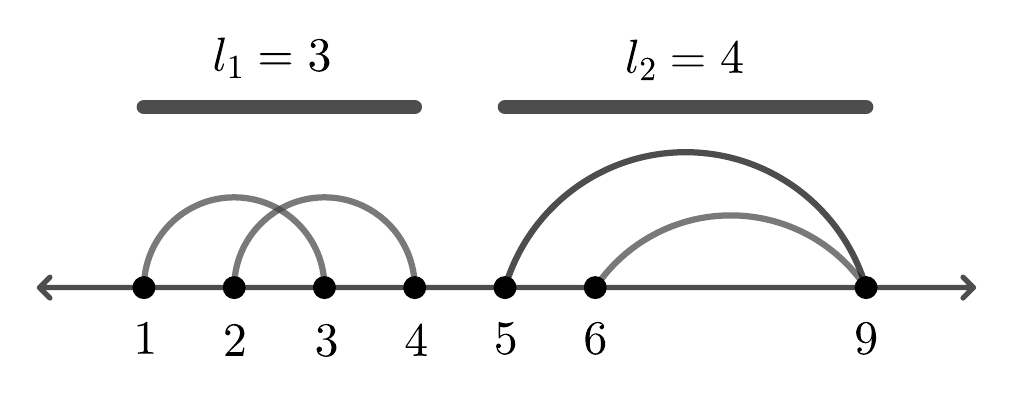}
\end{center}
\vspace{-.75cm}
\caption{\small The graph $G_t$ is embedded on the real line (self-loops not shown)
and $E_{s,t}= 3^s+4^s$. \label{fig-senergy}}
\end{figure}

Given $x(0)\in \mathbb{R}^n$, we write $x(t+1)= P_t x(t)$ for $t\geq 0$,
and refer to the $i$-th coordinate $x_i(t)$ as the position of agent $i$ at time $t$.
The vector $x(t)$ provides an embedding of the graph $G_t$ over the reals.
The union of the embedded edges of $G_t$ forms disjoint
intervals, called {\em blocks}.
Let $l_1(t),\ldots, l_{k_t}(t)$ be the lengths of these blocks
and put $E_{s,t}= \sum_{i=1}^{k_t} l_i(t)^s$, with $s\in (0,1]$.
For example, if, besides its self-loops, $G_t$ has four edges embedded 
as the intervals $[1, 3]$, $[2, 4]$, $[5, 9]$, and $[6,9]$, then there are $k_t=2$
blocks $[1, 4]$ and $[5,9]$; hence $E_{s,t} = 3^s + 4^s$
(Figure~\ref{fig-senergy}).

\vspace{0.4cm}

\vspace{-.25cm}
\vspace{0.5cm}

We define the {\em $s$-energy} $E_s= \sum_{t\geq 0} E_{s,t}$
and denote by $\mathcal{E}_{s}$ the supremum 
of $E_s\,$ over all $n$-vertex graph sequences $(G_t)_{t\geq 0}$
and any initial position  $x(0) \in [0,1]^n$.
If we limit the sequences to include only graphs with at most $m$ connected
components, then we denote the supremum by $\mathcal{E}_{m,s}$;
note that $\mathcal{E}_{s} = \mathcal{E}_{n-1,s}$.
The best known bound states that 
$\mathcal{E}_{s} \leq (3/\rho s)^{n-1}$~\cite{chazelle-Energ2-2019}.
We prove the following results, all of which assume 
the same fixed weight threshold~$\rho\in (0,1/2]$ and $s\in (0,1]$.
The superscripts below indicate the type of averaging systems considered.

\vspace{0.35cm}
\begin{theorem}\label{s-Energ-Bounds}
$\!\!\! \,\,$
\begin{enumerate}
\item[\rm{(a)}] \ \ 
$\mathcal{E}_{m,1} \leq 3en (1/\rho)^{\lfloor n/2\rfloor}$
and \
$\mathcal{E}_{m,s}\leq (c/s)^m(1/\rho)^{n-1}$,
where $c=O(mn)$.
\item[\rm{(b)}] \ \ 
$\mathcal{E}_{m,s}^{\text{\sc rev}}
\leq (cn^2/\rho s)^m$, for constant $c>0$.
\item[\rm{(c)}] \ \ 
$\mathcal{E}_{1,s}^{\, \text{\sc exp}}
= O(d^3/s h^2)$ \ and \ 
$\mathcal{E}_{m,s}^{\, \text{\sc exp}}
\leq (c/s)^m$, where $c= O(d^3mn/h^2)$.
\item[\rm{(d)}] \ \
$\mathbb{E}\, \mathcal{E}_s^{\, \text{\sc ran}}
\leq c/s$, for constant $c>0$.
\end{enumerate}
\end{theorem}
\smallskip
\vspace{0.35cm}

The bounds in (a) apply to any initial position $x(0) \in [0,1]^n$
or, equivalently, an initial diameter at most~1.
For notational convenience, we use a slightly different assumption for the cases (b, c, d):
we assume that the initial (scaled) variance $\sum_i q_i ( x_i(0)- \mu )^2$ is at most 1,
where $\mu= \sum_i \bar q_i x_i(0)$ is the mean position at time~$0$;
here, $\bar q= q/\|q\|_1$ is the stationary distribution of $P_t$.
To see why this different assumption can be made without loss of generality,
we note that unit diameter implies a variance of at most $n/\rho$
and, conversely, unit variance implies a diameter bounded by~2.

To our knowledge, the $s$-energy is the only tool
we have at this moment for proving convergence bounds for 
time-varying averaging systems.
It can also be used for random walks in {\em fixed} undirected graphs:
in fact, all known mixing times can be expressed (and rederived) in terms
of the $s$-energy.
The true utility of the concept, however, is in regard to {\em time-varying} graphs.
The proof of Theorem~\ref{s-Energ-Bounds} (given in~\S\ref{sec:proofs}) is quite novel.
It departs radically from typical proofs of mixing times. One of its original features
is its algorithmic nature. The stepwise contribution of an averaging system's orbit
to its $s$-energy is modeled by an algorithmic process much in the
spirit of amortized analysis.

\vspace{0.4cm}
\paragraph{Convergence rates.}

Let $T_\eps$ be the number of timesteps $t$ at which
$G_t$ has an embedded edge of length at least $\eps>0$. 
We denote by $T_{m,\eps}$ the maximum value of $T_\eps$
over all $n$-vertex sequences $(G_t)_{t\geq 0}$ with no more than
$m$ connected components.\footnote{Note that an adversary can always insert the identity matrix repeatedly and delay convergence at will, so counting all the steps prior to approaching the attractor
is not an option.}
We use superscripts as we did before.

\vspace{0.4cm}
\begin{theorem}\label{ConvergeUB}
$\!\!\! .\,\,$
For constant $c>0$ and any positive $\eps$ small enough,
\begin{equation*}
\begin{cases}
\ \  T_{m, \eps}\leq c (1/\rho)^{n-1}\big(mn\log \frac{1}{\eps}\big)^m 
\vspace{0.1cm} \\
\ \  T^{\text{\sc rev}}_{m, \eps}\leq c \big(\frac{n^2}{\rho}  \log \frac{1}{\eps} \big)^m 
\vspace{0.1cm} \\
\ \  T^{\text{\sc exp}}_{m, \eps}\leq  c \big( \frac{d^3mn}{h^2} \log \frac{1}{\eps}\big)^m 
\ \ and \ \ \ T^{\text{\sc exp}}_{1, \eps}\leq   \frac{cd^3}{h^2} \log \frac{1}{\eps}
 \vspace{0.1cm} \\
\ \  \mathbb{E}\, T^{\text{\sc ran}}_{1, \eps} \leq   c \log \frac{1}{\eps}
\end{cases}
\end{equation*}
\end{theorem}
\proof
By the definition of $T_{m,\eps}$,
the $s$-energy is bounded below by $\eps^s T_{m,\eps}$; 
hence, by Theorem~\ref{s-Energ-Bounds},
$$
T_{m,\eps}\leq \inf_{0<s\leq 1} \eps^{-s} \mathcal{E}_{m,s}
\leq  \inf_{0<s\leq 1} \eps^{-s} (c/s)^m(1/\rho)^{n-1} ,
$$
where $c=O(mn)$. Plugging in $s= b/\log(1/\eps)$, for constant $b>0$,
proves the claimed bound
on $T_{m,\eps}$. The three other upper bounds are derived similarly.
\hfill $\Box$
\proofend

Note that the case $m=1$ in $T_{m, \eps}$ and $T^{\text{\sc rev}}_{m, \eps}$
recovers the classic mixing times for Markov chains ($P_t=P$) and, in particular,
the polynomial vs. exponential gap between general and reversible chains.
It also recovers the usual logarithmic bound for expanders of constant degree~\cite{LPW}.

\section{Energy Bounds for Averaging Systems}\label{sec:proofs}

In this section, we prove Theorem~\ref{s-Energ-Bounds}. The four parts correspond to different
structural assumptions on the averaging system. Part~(a) is the most general case and contains
the core argument of the paper. Parts~(b)-(d) introduce additional structure that 
allows alternative techniques and leads to sharper bounds.

\subsection{General averaging systems: Theorem~\ref{s-Energ-Bounds} (a) }\label{thma}

An averaging system is a special case of a 
{\em twist system}~\cite{chazelle-Energ2-2019}.
The latter is easier to analyze so we turn our attention to~it.
As opposed to an averaging system, where the next step
is uniquely specified, a twist system is nondeterministic.
Instead of a fixed rule spelling out the chronological transitions,
we are given constraints at each timestep, and any transition
that conforms to them is valid. 

Relabel the agents so their
positions $x_1\leq \cdots \leq x_n$ appear in sorted order at time~$t$.
A twist system moves them to positions
$y_1 \leq \cdots \leq y_n$ at time $t+1$ in such a way that
\begin{equation}\label{twist-cond}
(1-\rho)x_u + \rho x_{\min\{i+1, v\}}\leq y_i\leq 
\rho x_{\max\{i-1,u\}} + (1-\rho)x_v,
\end{equation}
for any $i$ in $[u,v]$ and $y_i=x_i$ otherwise.
We repeat this step indefinitely. 
Twist systems are highly nondeterministic.
At each step, a new interval $[u,v]\subseteq [n]$, called a {\em block}, is picked
and the agents' motion is only constrained
by~(\ref{twist-cond}) and the need to maintain their ranks (ie, agents never cross).
It is not entirely obvious that twist systems should always exist:
in other words, that the constraints imposed in~(\ref{twist-cond}) are always satisfiable.
This is a direct consequence of Lemma~\ref{SAS<Twist}, which we state and prove below.

For the purposes of this work,
we extend the concept to {\em $m$-twist systems} by stipulating, 
at each time $t$, a partition of $[n]$ into up to $m_t$ blocks
$[u_{t,l}, v_{t,l}]$ ($1\leq l\leq m_t\leq m$).
Each agent is now subject to~(\ref{twist-cond}) within its own enclosing block.
We define the $s$-energy $E_s^{\text{\sc tw}}$
of a twist system as we did with an averaging system by adding
together the $s$-th powers of all the block lengths. We use the same notation
with the addition of the superscript {\sc tw}.

\begin{lemma}\label{SAS<Twist}
$\!\!\! .\,\,$
An averaging system with at most $m$ connected components at any time
can be interpreted as an $m$-twist system with the same $s$-energy. 
\end{lemma}
\proof
Fix an averaging system and let $(x_i)_{i=1}^n$ and $(y_i)_{i=1}^n$ be
the positions of the agents at times~$t$ and $t+1$, given in nondecreasing order.
We denote by $x_i'$ the position of agent $i$ at time $t+1$.
Let $[x_u,x_v]$ be a block at time $t$
(ie, an interval of the union of the embedded edges of $G_t$).
Pick $k<v$ and write $z= \rho x_k + (1-\rho)x_v$.
All the diagonal elements of $P_t$ are at least $\rho$; hence
$x'_i\leq z$, for all $i\leq k$, and $y_k\leq z$.  In fact,
we even have $y_{k+1} \leq z$. 
Indeed, the embedded edges of $G_t$ cover all of $[x_u,x_v]$, so at least one of
them, call it $(l,r)$, must join $[u,k]$ to $[k+1,v]$; hence $x'_r\leq z$. Our claim follows.
This proves that, for all $i\in (u,v]$, $y_i\leq \rho x_{\max\{i-1,u\}} + (1-\rho)x_v$.
We omit the case $i=u$ and the mirror-image inequality, which repeat the same argument.
Summing up all the powers $(x_v-x_u)^s$ shows the equivalence between the two $s$-energies.
\hfill $\Box$
\proofend

We may assume that the agents stay within $[0,1]$.
We begin with the bound $\mathcal{E}_{m,1} \leq 
\mathcal{E}_{m,1}^{\text{\sc tw}} \leq 3en (1/\rho)^{\lfloor n/2\rfloor}$.

\medskip
\begin{center}
\fbox{\rm Case $s=1$}
\end{center}
\bigskip
We prove a stronger result by bounding 
$K_{t}(z) \!:= \sum_{k=1}^n \left(x_{v(k)}- x_k\right) z^k$, where
$v(k)= v_{t,l}$ for $l$ such that $k\in [u_{t,l}, v_{t,l}]$.
As usual, $0\leq x_1\leq\cdots\leq x_n\leq 1$ denotes the sorted positions
of the agents at time~$t$; we omit $t$ for convenience but it is understood throughout.
We define the {\em weighted $1$-energy} 
$K(z)  = \sum_{t \geq 0} K_{t}(z)$
and, finally, $\mathcal{K} (z) = \sup K(z)$.
As long as $z\geq 1$, the $1$-energy is obviously dominated
by its weighted version. We improve this crude bound via a symmetry argument:

\begin{lemma}\label{weighted-E}
$\!\!\! .\,\,$
For any $z\geq 1$,  
$\mathcal{E}_{m,1}^{\text{\sc tw}}
\leq 2 z^{-\nu}\, \mathcal{K}(z)$, where $\nu= \lceil n/2 \rceil$.
\end{lemma}
\proof
We define the mirror image of $K_{t}$ as
$\bar K_{t}(z)= \sum_{k=1}^n \left(x_k- x_{u(k)}\right)  z^{n-k+1}$, where
$u(k)$ is the left counterpart of $v(k)$. We have

\begin{equation*}
\begin{split}
E_{1,t}^{\text{\sc tw}}
&\leq  \sum_{k\leq \nu}  \left(x_k- x_{u(k)}\right)
     + \sum_{k\geq \nu} \left(x_{v(k)}- x_k\right)  \\
&\leq  z^{-\nu} \sum_{k\leq \nu}  \left(x_k- x_{u(k)}\right) z^{n-k+1}
     +  z^{-\nu}  \sum_{k\geq \nu} \left(x_{v(k)}- x_k\right) z^k \\
&\leq z^{-\nu} \bigl(  \bar K_{t}(z) + K_{t}(z) \bigr).
\end{split}
\end{equation*}
Because $\mathcal{K}(z) = \sup K(z)$, the lemma then follows by summing up 
over all $t \geq 0$.
\hfill $\Box$
\proofend

We define the polynomial\footnote{Not to be confused with the stochastic matrix $P_t$
used earlier.}
$P_t(z)= \sum_{k=1}^n x_k z^k$ for $z> 1/\rho$
and exploit two simple but surprising facts:
$P_t(z)$ cannot increase over time;\footnote{Recall that $x_k$ depends on $t$.
Note also that, among the $n$ agents,
rightward motion within $[0,1]$ might greatly outweigh the leftward kind.
Thus, if most of the $x_i$'s keep growing, how can $P_t(z)$ not follow suit?
The point is that $P_t(z)$
puts weights exponentially growing on the right, so their leftward motion,
outweighed as it might be, will always dominate with respect to $P_t(z)$.
This balancing act between left and right motion is 
the core principle of twist systems.}
and, at each step, the drop from $P_t(z)$ to $P_{t+1}(z)$ is at least 
proportional to $K_{t}(z)$.
Thus, we develop a discrete version of the inference: $dP_t/dt \leq - c K_{t}$ implies
$$
 \int_{t\geq 0} c K_{t} \leq  - \int_{t\geq 0} \frac{dP_t}{dt}
\leq P_1\, .
$$ 
\begin{lemma}\label{PQ}
$\!\!\! .\,\,$
For any $z>1/\rho$, $P_t(z)- P_{t+1}(z)\geq (\rho z -1) K_{t}(z)$.
\end{lemma}
\proof
The inequality is additive in the number of blocks so we can assume
there is a single block $[u,v]$ at time~$t$.
Using the notation of~(\ref{twist-cond}), we have 
$y_k\leq \rho x_{\max\{k-1,u\}} + (1-\rho)x_v$; hence
\begin{equation*}
\begin{split}
P_t(z) - P_{t+1}(z)
&= \sum_{k=u}^v (x_k-y_k)z^k
\geq  \sum_{k=u}^v \left(x_k-\rho x_{\max\{k-1,u\}} - (1-\rho)x_v \right)z^k \\
&\geq  (\rho-1)(x_v- x_u)z^u + \sum_{k=u+1}^v \rho  (x_v- x_{k-1})z^k
-  \sum_{k=u+1}^v (x_v-x_k) z^k \\
&\geq  (\rho-1)(x_v- x_u)z^u + \sum_{k=u}^{v-1} \rho z (x_v- x_k)z^k
-  \sum_{k=u+1}^v (x_v-x_k) z^k \\
&\geq  (\rho-1)(x_v- x_u)z^u + \sum_{k=u}^v (\rho z-1) (x_v- x_k)z^k
 + (x_v- x_u)z^u \\
&\geq  (\rho z -1) K_t(z) + \rho (x_v- x_u)z^u.
\end{split}
\end{equation*} 
\hfill $\Box$
\proofend

The lemma implies that 
\begin{equation}\label{Qz-ub}
(\rho z -1) K(z)
=   \sum_{t \geq 0} (\rho z -1) K_{t}(z) \leq P_0(z)
\leq \sum_{k=1}^n z^k = \frac{z^{n+1}-z}{z-1} \, .
\end{equation}
With $z= (1+\eps)/\rho$, $\eps= 1/(n-\nu+1)$, and $n>2$, we find that
$$
z^{-\nu}K(z)
\leq \frac{z^{n}-1}{(z-1)(\rho z -1)z^{\nu-1}}
\leq 2\rho^{\nu -n}e^{(n-\nu+1)\eps}/\eps 
\leq 2e(n-\nu+1) \rho^{\nu -n}\leq \frac{3en}{2} \rho^{\nu-n}.
$$
The case $s=1$ of Theorem~\ref{s-Energ-Bounds}~(a) follows immediately
from Lemma~\ref{weighted-E}.  Finally, for $n=2$, we verify that
$\mathcal{E}_{m,1}^{\text{\sc tw}} = \sum_{k\geq 0} (1-2\rho)^k= 1/2\rho$.
\hfill $\Box$
\proofend

\begin{center}
\fbox{\rm Case $s<1$}
\end{center}
\medskip

The previous argument relied crucially on the linearity of the $1$-energy.
If $s<1$,  the $s$-energy gives
more relative weight to small lengths, so we need
a different strategy to keep the scales separated.
We omit the superscript {\sc tw} below but it is understood.
We use a threshold $\delta$ which, though set to $1/3$, is best kept
as $\delta$ in the notation. 

\bigskip
\noindent
{\em A recurrence relation.}\ \ \ 
Let $T_\delta$ be the number of steps at which
the diameter remains above $1-\delta$; note that 
these steps are consecutive and $T_\delta$ might be infinite.
By scaling, we find that
$E_s \leq F_{\!s} + (1-\delta)^s E_s$,
where $F_{\!s} = \sum_{t\leq T_\delta} E_{s,t}$. 
Since $(1-\delta)^s \leq 1 - \delta s$, for $\delta,s\in [0,1]$, we have
\begin{equation}\label{BoundEF}
E_s \leq  (\delta s)^{-1}  F_{\!s}.
\end{equation}

\noindent
If $m=1$, then $ (1-\delta)F_{\!s} \leq  (1-\delta)T_\delta \leq F_1\leq \mathcal{E}_{1,1}$.
Thus, by~(\ref{BoundEF}) and the previous section,
\begin{equation}\label{Bound:m=1}
\mathcal{E}_{1,s} \leq \frac{3en}{\delta(1-\delta)}(1/s)(1/\rho)^{\lfloor n/2\rfloor}.
\end{equation}

\medskip
\noindent
We consider the case $m>1$. 
Fix $t\leq T_\delta$; if $m_t>1$, let $j$ maximize $x_{j+1}(t) - x_j(t)$ over all
$j= v_i$ and $i< m_t$ (break ties by taking the smallest $j$).
This corresponds
to the maximum distance between consecutive blocks.
For this reason, we call $(j,j+1)$ the {\em max-gap} at time~$t$.
We say that $t$ is {\em ungapped} if $m_t=1$ 
or $x_{j+1}(t) - x_{j}(t)\leq \delta/m$;
it is {\em gapped} otherwise.
Assuming that $t$ is gapped, 
let $(j, j+1)$ be its max-gap and write
$\zeta_t = \min_{k} \bigl\{\, t< k \leq T_\delta \,|\,  \exists\, l  : u_{k,l}\leq j < v_{k,l} \, \bigr\}$:
If the set is empty, we set $\zeta_t =T_\delta$; else $l$ is unique and we denote it by $l_t$.
We call the interval $[t,\zeta_t]$ a {\em span} and the block $l_t$,
if it exists, its {\em cap}.
We note that $x_{j+1}(k)- x_j(k)$ cannot decrease during the times
$k=t,\ldots, \zeta_t$. This shows that a cap covers a length greater than $\delta/m$.

We begin with a few words of intuition. The energetic contribution of 
an ungapped time $t$ is easy to account for:  It is at most $m$. On the
other hand, the $1$-energy is at least the diameter minus the added length of the gaps between blocks,
which amounts to at least $1-2\delta\geq 1/3$;
in other words, $E_{s,t}\leq 3m E_{1,t}$. Summing up over all ungapped 
times and plugging in our bound for $s=1$ gives us the desired result.
Accounting for gapped times is more difficult,
as it requires dealing with small scales. If we had only one span,
we could simply split the system into two decoupled subsystems 
and set up a recurrence relation. 
The problem is that the presence of $k$ capped spans would force us to repeat the
recursion $k-1$ times. With no apriori bound on $k$, this approach is not too promising.
Instead, we make a bold move: We argue that, because a cap is longer than~$\delta/m$,
its own 1-energy contribution (ie, its length) is large enough
to ``pay'' for the $s$-energy of its entire span. This is not quite right, of course,
but one can fix the argument by using the weighted 1-energy of the cap
and upscaling it suitably.
Once again, this reduces the problem to the case $s=1$, so 
our method is, in effect, a linearization.
Here is the proof.

We partition the times between 0 and $T_\delta$ into two subsets $\mathcal{G}$
and $\mathcal{U}  \!:= [0,T_\delta]\setminus \mathcal{G}$, each one supplied
with its own energetic accounting scheme.
We form $\mathcal{G}$ by greedily extracting a maximal set of nonoverlapping spans
and taking their union.

\vspace{0.3cm}

{\small
\par\medskip
\renewcommand{\sboxsep}{0.5cm}
\renewcommand{\sdim}{0.8\fboxsep}
\begin{center}
\shabox{\parbox{8cm}{
\begin{itemize}
\item[\text{1.}]
\hspace{0.1cm}
$\mathcal{G} \leftarrow \emptyset$ and $t' \leftarrow 0$;
\item[\text{2.}]
\hspace{0.1cm}
{\bf if } $t\leftarrow \min \, \bigl\{\, \text{gapped } i \, |\, t' \leq i \leq T_\delta \, \bigr\}$ exists 
\item[\text{3.}]
\hspace{0.5cm}
{\bf then } $\mathcal{G}\leftarrow \mathcal{G} \cup \bigl\{\,i\,|\, t\leq i\leq  \zeta_t\, \bigr\}$;
\item[\text{4.}]
\hspace{1.35cm}
{\bf if } $\zeta_t<T_\delta$ {\bf then } $t'= \zeta_t +1$; go to 2;
\end{itemize}
}}
\end{center}
\par
}
\vspace{.5cm}

\noindent
We postulate that, for any $0<s\leq 1$, and any number of agents
$j\leq n$,
\begin{equation}\label{Rstmnz}
\mathcal{E}_{m,s} \leq c_{m}(1/s)^m( 1/\rho)^{j-1} \, ,
\end{equation}
and we derive a recurrence relation for $c_{m}$ (for given $n$).

\begin{itemize}
\item
{\em Accounting for $\mathcal{G}$:}
In line 3, let $k=\zeta_t$ and $(j, j+1)$ be the corresponding max-gap.
Suppose that the span $[t,k]$ is capped.
The absence of an interval including $j$ and $j+1$ during 
$[t, k -1]$ implies that
$\sum_{t\leq l\leq k} E_{s,l} \leq L+R +m$, where $L$ and $R$
denote the $s$-energy of systems with at most $m-1$ connected components.
For reasons we address below, we may assume that $L$ is dominant;
hence $R\leq L\leq \mathcal{E}_{m-1,s}$.\footnote{The inequality 
relies on the (easy) fact that the maximum $s$-energy
grows monotonically with the number of agents. This is not even
needed, however, if we redefine $\mathcal{E}_{m,s}$ as the maximum
$s$-energy over all systems with {\em at most} $n$ agents
and then reason with the value $n'\leq n$ that achieves the maximum.}
It follows that
\begin{equation}\label{Et-k}
\sum_{l=t}^k  E_{s,l} \leq   2\mathcal{E}_{m-1,s} +m
\leq 3c_{m-1}(1/s)^{m-1}( 1/\rho)^{j-1}.
\end{equation}
Note that we (safely) assume $c_{m-1}\geq m$.
Using the shorthand $v$ for $v_{k,l_t}$, we have
$x_v(k)- x_j(k) \geq x_{j+1}(k) - x_j(k) 
\geq x_{j+1}(t) - x_j(t) > \delta/m$.
We add the artificial multiplier $x_v(k)- x_j(k)$ to~(\ref{Et-k})
to make the right-hand side resemble $K(z)$.
Recall that $\delta=1/3$; assuming that $z>1/\rho$ 
from now on, we have

\begin{equation}\label{Gap-E}
\sum_{l=t}^k  E_{s,l}
\leq  B\bigl(x_v(k) - x_j(k)\bigr) z^j,
\hspace{0.4cm} \text{\small{with}} \hspace{0.2cm}
B=  9 c_{m-1} m\rho (1/s)^{m-1}
\end{equation}

The set $\mathcal{G}$ is a union of spans. If $\zeta_t= T_\delta$, the last span might not be capped.
If so, remove it from $\mathcal{G}$ and call the resulting set $\mathcal{G}'$.
Summing up, we find that
$\sum_{t\in \mathcal{G}'} E_{s,t} \leq  B \sum_{t\in \mathcal{G}'} K_t(z)$.
If the last span is uncapped then $\zeta_t=T_\delta$ and no block
contains both $j$ and $j+1$ in the span $[t, T_\delta]$.
The $s$-energy expended in that span is 
thus of the form $L+R\leq 2\mathcal{E}_{m-1,s}$.

\item
{\em Accounting for $\mathcal{U}$:}
Only ungapped times belong to $\mathcal{U}$, so the 1-energy at time $t\in \mathcal{U}$
is at least $1-\delta - (m_t-1)\delta/m\geq 1/3$.
On the other hand,
$E_{s,t} \leq m\leq 3m E_{1,t}
\leq 3m \rho K_t(z) \leq B K_t(z)$.
\end{itemize}

\noindent
Set $z= (1+\eps)/\rho$, for $\eps>0$.
Putting all of our bounds together, we have  
\begin{equation*}
\begin{split}
F_{\!s}
= \sum_{t\leq T_\delta} E_{s,t}
&\leq B\sum_{t\in \mathcal{G}'} K_t(z) +
2\mathcal{E}_{m-1,s}
+  B\sum_{t\in \mathcal{U}} K_t(z) \\
&\leq 2 c_{m-1} \,  (1/s)^{m-1} (1/\rho)^{n-1}
+ BK(z) .
\end{split}
\end{equation*}
Actually, the exponent to $1/\rho$ can be reduced to $n-2$, but this is immaterial.
By~(\ref{Qz-ub}),
$$B K(z)\leq  
18 c_{m-1}  (1/s)^{m-1} (1/\rho)^{n-1} e^{(n+1)\eps} (m/\eps).
$$
Setting $\eps= 1/(n+1)$ gives us, for some constant $d>0$,
\begin{equation*}
F_{\!s}
\leq  c_{m-1} (dmn) (1/s)^{m-1} (1/\rho)^{n-1}  .
\end{equation*}

We tie up the loose ends by arguing that it was legitimate to assume that
$L\geq R$.
The point is that individual values of $L$ and $R$ do not matter: only their sums do.
Thus, if the $R$s outweigh the $L$s,
we restore the dominance of the $L$s by flipping the system around.
Finally, by~(\ref{BoundEF}),
$E_s\leq (3/s) F_{\!s}$; and so, by~(\ref{Bound:m=1}),
Theorem~\ref{s-Energ-Bounds}~(a)
follows from the recurrence:
$c_{1}= O(n)$ and $c_{m}\leq 3dmn c_{m-1}$, for $m>1$.
\hfill $\Box$
\proofend

\paragraph{Lower bounds for twist systems.}

We begin with the case $m=1$.
Assume that $n=2k+1$.
At time $t=1$, we have $x_k=- x_{-k}= 1/2$
and $x_i=  -x_{-i} = \frac{1}{2}(1-\rho^i) $, for $0\leq i<k$
and $\rho$ small enough. 
For $t>1$, we set $x_i(t)= (1-\rho^k)x_i(t-1)$.
The agents are labeled $-k,\ldots, k$ from left to right.
It is easily verified that this constitutes a twist system for the block $[-k,k]$
with initial unit diameter.
The $s$-energy $E$ 
is $\bigl(1-\rho^k\bigr)^s E +1$ so, for constant $c>0$,
\begin{equation}\label{E-tw-m=1}
E \geq (c/s)(1/\rho)^{\lfloor n/2\rfloor}.
\end{equation}
If $n=2k$, we set $x_k=- x_{-k}= 1/2$
and $x_i=  -x_{-i} = \frac{1}{2}(1- 2\rho^i) $, for $1\leq i<k$.
For $t>1$, we set $x_i(t)= (1-2\rho^k)x_i(t-1)$ and 
rederive~(\ref{E-tw-m=1}).

For the general case, we describe the evolution of an $m$-block twist system
with $n$ agents, and denote its $s$-energy by $F(n,m)$:
It is assumed that $n-1$ agents are positioned at $0$ at time~1
and the last one is at position~$1$.
If $m=1$, we apply the previous construction after shifting the initial
interval from $[-0.5, 0.5]$ to $[0,1]$.
The initial positions still do not match,
but we note that, in a single step, we can move the agents anywhere we want in
the interval $[\rho, 1-\rho]$ while respecting the constraints of a twist system.
This gives us
$F(n,1)= 1 + (1-2\rho)^s E$.
By adjusting the constant $c$ in~(\ref{E-tw-m=1}), the same lower bound still holds.

For $m>1$, at time $1$, we move the agents $n-1$ and $n$ to positions
$\rho$ and $1-\rho$, respectively, and we leave the others (if any) at position 0.
We then use an $(m-1)$-block twist system recursively 
for the agents $1,\ldots, n-1$. This brings these agents to
a common position\footnote{To keep the time finite,
we can always force completion in a single step
once the agents are sufficiently close to each other.} in $[0,\rho]$.
This gives us the recurrence relation:
$F(n,m)\geq 1+ \rho^s F(n-1, m-1) + (1-2\rho)^sF(n,m)$; hence, by induction,
for constant $c>0$,
\begin{equation}\label{E-tw-all-m}
F(n,m)\geq (c/s)^m \rho^{(m-1)s} (1/\rho)^{\lfloor (n+m-1)/2\rfloor}.
\end{equation}
The $s$-energy is often used for small $s$, so we state 
the case of $s=O\bigl(1/m\log \frac{1}{\rho}\bigr)$,
which matches the bound
of Theorem~\ref{s-Energ-Bounds}~(a) for $m=n-1$.

\begin{theorem}\label{s-energy-TW-LB}
$\!\!\! .\,\,$
$\mathcal{E}_{m,s}^{\text{\sc tw}} \geq
     (c/s)^m (1/\rho)^{\lfloor (n+m-1)/2\rfloor}$,
for constant $c>0$, small enough $\rho$ and $s= O\bigl(1/m\log \frac{1}{\rho}\bigr)$.
\end{theorem}

\subsection{Reversible averaging systems: Theorem~\ref{s-Energ-Bounds}~(b)}\label{RAS}

Recall that, in a reversible averaging system {\sc rev},
$P_t= \text{diag}(q)^{-1} M_t$, where
$M_t$ is symmetric with nonzero entries at least~1
and $q= M_t\mathbf{1}\preceq  \mathbf{1}/\rho$.
The vector $q$ is proportional to the stationary distribution
of the Markov chain induced by $P_t$. By reversibility,
we have $q_i (P_t)_{ij}= q_j (P_t)_{ji}$.
Write 
$\langle x, y \rangle_q : = \sum_i q_i x_i y_i$
and $\|x\|_q^2:= \langle x, x \rangle_q$.
The (scaled) variance defined in~\S\ref{sec:intro} can be rewritten
as $\|x-  \hat x\, \|_q^2$ where
$\hat x= \|q\|_1^{-1}\langle x, \mathbf{1}\rangle_q\mathbf{1}$
and $x$ is shorthand for $x(0)$.
The convergence rate of the attracting dynamics
is captured by a variant of the Dirichlet form:
\begin{equation}\label{DirichletDef}
D_t= \sum_{i}\max_{j:\, (i,j)\in G_t} \bigl( x_i(t)- x_j(t) \bigr)^2.
\end{equation}
We omit the index $t$ below for clarity.

\begin{lemma}\label{dirich}
$\!\!\! .\,\,$
\ 
$\|Px\|_q^2  \leq \|x\|_q^2 - D/2\,$, for any $x\in \mathbb{R}^n$.
\end{lemma}
\proof
Write $\delta_{ij}= x_i-x_j$ and $\mu_i = \sum_{j} p_{ij}  \, \delta_{ij}$.
Fix $i$ and pick any $k$ such that $m_{ik}>0$.
By Cauchy-Schwarz and $m_{ii}\geq 1$ (because nonzero entries of $M$ are at least $1$), we have
\begin{align*}
\delta_{ik}^2 &= \bigl( (\mu_i- \delta_{ii}) + (\delta_{ik}-\mu_i) \bigr)^2
\leq 2(\delta_{ii}-\mu_i)^2 + 2(\delta_{ik}-\mu_i)^2
   \leq 2\sum_{j} m_{ij}  \, (\delta_{ij}- \mu_i)^2;
\end{align*}
hence,
\begin{align*}
\|x\|_q^2- \|Px\|_q^2  
&=  \sum_i q_i x_i^2 -  \sum_i q_i \,\biggl(x_i+ \sum_j p_{ij}\delta_{ji}\biggr)^2
= - \sum_i q_i \, \biggl( 2 x_i \sum_j p_{ij}\delta_{ji} +  \mu_i^2  \biggr) \\
& =  \sum_i q_i \, \biggl( \, \sum_{j} p_{ij}  \, \delta_{ij}^2   -  \mu_i^2 \, \biggr) 
=  \sum_{i,j} m_{ij}  \, (\delta_{ij}- \mu_i)^2 
\geq \frac{1}{2}\sum_{i}  \max_{j: m_{ij}>0} \delta_{ij}^2,
\end{align*}
with the last equality expressing the identity for the variance:
$\mathbb{E} X^2 - (\mathbb{E}\, X)^2
= \mathbb{E}[X- \mathbb{E}\,X]^2$.
\hfill $\Box$
\proofend

Write $G_{\leq t}$ as the union of all the edges in $G_0,\ldots, G_t$,
and let $t_c$ be the maximum value of~$t$ such that
$G_{\leq t}$ has fewer connected components than $G_{\leq t-1}$;
if no such $t$ , set $t_c=0$.

\begin{lemma}\label{cover-length}
$\!\!\! .\,\,$
If $G_{\leq t_c}$ is connected, then
$\sum_{t\leq t_c} D_t \geq \, \rho n^{-2}\,  \|x-\hat x\, \|_q^2\,$.
\end{lemma}
\proof
Let $G_{t_0}$ denote the graph over $n$ vertices with no edges.
We define $t_1,\ldots, t_c$ as the sequence of times $t$ at which
the addition of $G_{t}$ reduces the number of connected components in $G_{\leq t-1}$.
At any time $t_k$ ($k>0$), 
the drop $d_k$ in the number of components can be achieved by $d_k$
edges from $G_{t_k}$.
Let $F_k$ denote such a set of edges: we can always
order $F_k$ so that every edge in the sequence
contains at least one vertex not encountered yet.
This shows that the sum of the squared lengths
of the edges in $F_k$ does not exceed $D_{t_k}$.
We note that 
$F:= F_1\cup\cdots\cup F_c$ forms a collection of $n-1$ edges
from (the connected graph) 
$G_{\leq t_c}$ and $F$ spans all $n$ vertices.

Consider the intervals formed by the edges in $F_k$ at time $t_k$,
for all $k\in [c]$. 
The union of these intervals covers the smallest interval $[a,b]$
enclosing all the vertices at time $0$ (and hence at all times). 
To see why, pick any $z$ such that $a<z<b$ and denote by $L$
and $R$ the vertices on both sides of $z$ at time $0$.   Neither set is empty and,
by convexity, both of them remain on their respective side
of $z$ until an edge of some $G_t$ joins $L$ to $R$.
When that happens (which it must since $G_{\leq t_c}$ is connected), the joining edge(s)
reduce(s) the number of components of $G_{\leq t-1}$ 
by at least one, so $F$ must grab at least one of them,
which proves our claim. Let $l_1,\ldots, l_{n-1}$ denote the
lengths of the edges of $F$ (at the time of their insertion).
By Cauchy-Schwarz, 
$\sum_{t\leq t_c} D_t \geq \sum_{i=1}^{n-1} l_i^2 \geq (b-a)^2/(n-1)$.
The lemma follows from the inequalities
$\|x-\hat x\, \|_q^2\leq \|q\|_1(b-a)^2\leq (n/\rho) (b-a)^2$.
\hfill $\Box$
\proofend

\smallskip

Assume that $G_{\leq t_c}$ is connected.
By Lemma~\ref{cover-length} and the telescoping use of Lemma~\ref{dirich},
\begin{equation}\label{telescope}
 \|x\|_q^2 -  \|x(t_c+1)\|_q^2
\geq \frac{1}{2} \sum_{t=0}^{t_c} D_t
\geq \frac{\rho}{2n^2} \, \|x-\hat x\, \|_q^2 \, .
\end{equation}

\noindent
Let $U(n,m)$ be the maximum $s$-energy of a {\sc rev}
with at most $n$ vertices and $m$ connected components at any time,
subject to the initial condition $\|x - \hat x\, \|_q^2 \leq 1$.
By shifting the system if need be, we can always assume that $\hat x= \mathbf{0}$.
By~(\ref{telescope}),
$\|x(t)\|_q^2$ shrinks by at least a factor of $\alpha:= 1- \rho/2n^2$ by time $t_c+1$.
A simple scaling argument shows that
the $s$-energy expanded after $t_c$ is at most $\alpha^{s/2} U(n,m)$.
While $t<t_c$ (or if $G_{\leq t_c}$ is not connected),
the system can be decoupled into two {\sc rev} systems,
each one with fewer than $m$ components.\footnote{Note that
each subsystem satisfies the required
inequalities about the $Q$ and $M$ entries; also, shifting
each subsystem so that $\hat x=0$ cannot increase
$\|x\|_q^2$, so its value remains at most 1.}
Since $\|x\|_q^2 \leq 1$, we have $|x_i|\leq 1$ for each $i$ and 
the diameter at any time is at most 2; therefore
$U(n,m) \leq \alpha^{s/2} U(n,m) + 2 U(n,m-1) + m2^s$.
It follows that 
\begin{equation}\label{U-recur}
U(n,m) \leq \frac{2}{1-\alpha^{s/2}}\, \bigl( U(n,m-1) + m \, \bigr);
\end{equation}
hence $U(n,m) = O(n^2/\rho s)^m$,
and Theorem~\ref{s-Energ-Bounds}~(b) follows.
\hfill $\Box$
\proofend

\vspace{1cm}

\paragraph{A lower bound for reversible systems.}

We begin with the case $m=1$.
The path graph $G$ over $n$ vertices has an edge $(i,i+1)$ for all $i<n$.
The Laplacian $L$ is $\text{diag} (u) -A$, where $A$ is the adjacency matrix of $G$
and $u$ is the degree vector $(1,2,\ldots, 2,1)$. We consider 
the reversible averaging system formed by the matrix $P_t= P= I - \rho L$.
By well-known spectral results on graphs~\cite{spielman19},
$P_t$ has a full set of $n$ orthogonal eigenvectors
$v_k$, where $v_k(i)=  \cos \frac{(i-1/2)k\pi}{n}$ for $i\in [n]$,
with its associated eigenvalues
$\lambda_k= 1- 2\rho\bigl(1-  \cos \frac{ k\pi}{n}\bigr) $, for $0\leq k<n$.
We require $\rho<1/4$ to ensure that $P$ is positive semidefinite.
We initialize the system with $x= (1,0,\ldots, 0)$ and observe that
the agents always keep their initial rank order, so the diameter $\Delta_{t}$ at time $t$
is equal to $(1,0,\ldots, 0,-1) P^{t} x$.
We verify that $\|v_1\|^2= n/2$. By the spectral identity 
$P^j= \sum_{k<n} \lambda_k^j  v_k v_k^T/  \|v_k\|^2$, 
we find that, for $t>0$ and $n>1$,
\begin{equation*}
\begin{split}
\Delta_t
&= \sum_{k=0}^{n-1} \lambda_k^{t}  v_k(1)\frac{v_k(1)-v_k(n)}{\|v_k\|^2}
= \sum_{{\rm\small odd}\,\,\, k}
\frac{2\lambda_k^{t}}{\|v_k\|^2} \Bigl(\cos\frac{k\pi}{2n} \Bigr)^2
\geq  \frac{2\lambda_1^{t}}{\|v_1\|^2} \Bigl(\cos\frac{\pi}{2n} \Bigr)^2
\geq \frac{2}{n}\lambda_1^{t}.
\end{split}
\end{equation*}
The $s$-energy is equal 
to $\sum_t \Delta_t^s\geq  (2/n)^s/(1-\lambda_1^s)
\geq bn^{2-s}/\rho s$, for constant $b>0$.

For the general case, we denote by $F(n, m)$ the $s$-energy of the 
system with initial diameter equal to 1. We showed that
$F(n,1)\geq  bn^{2-s}/\rho s$.
We now describe the steps of the dynamics for $m>1$.
To simplify the notation, we assume that
$\nu:= n/m$ is an integer.\footnote{This can be relaxed with a simple
padding argument we may omit.}
For $i\in [m]$, let $C_i$ be the path linking vertices $[(i-1)\nu +1, i\nu]$.

\begin{enumerate}
\item
At time $t=1$, the vertices of $C_1$ are placed at position $0$
while all the others are stationed at $1$. 
The paths $C_1$ and $C_2$ are linked together into a single path
so the system has $m-1$ components. Vertices $\nu$ and $\nu+1$ move
to positions $\rho$ and $1-\rho$ respectively while the others do not move at all.
The $s$-energy expended during that step is equal to 1.
\item
The system now consists of the $m$ paths $C_i$.  We apply the case $m=1$
to $C_1$ and $C_2$ in parallel, which expends $s$-energy
equal to $2\rho^s F(\nu, 1)$. All other vertices stay in place.
The transformation keeps the mass center invariant, so 
the vertices of $C_1$ and $C_2$ end up at positions $\rho/\nu$
and $1- \rho/\nu$, respectively.\footnote{To keep the time finite,
we can always force completion in a single step
once the agents are sufficiently close to each other
and use a limiting argument.} 
\item
We move the vertices in $C_i$ for $i\geq 2$ by applying
the same construction recursively for fewer than $m$ components.
The vertices of $C_1$ stay in place.
The $s$-energy used in the process is equal to 
$(\rho/\nu)^s F(n-\nu, m-1)$ 
and the vertices of $C_2,\ldots, C_m$
end up at clustered at position $1- \frac{\rho}{n-\nu}$.
\item
We apply the construction recursively to the $n$ vertices,
which uses up a quantity of $s$-energy equal to 
$\bigl( 1-\frac{\rho}{\nu} - \frac{\rho}{n-\nu} \bigr)^s F(n,m)$.
\end{enumerate}
Putting all the energetic contributions together, we find that,
for constants $b', c>0$,
\begin{equation*}
\begin{split}
F(n,m) &\geq 
1 +
\frac{2b \nu^{2-s}}{s \rho^{1-s}}  + 
\left( \frac{\rho}{\nu} \right)^s F(n-\nu, m-1)
+\left( 1-\frac{2\rho}{\nu}  \right)^s F(n,m) \\
&\geq \frac{b'}{s\rho^{1-s}}\left(\frac{n}{m}\right)^{1-s}F(n-\nu, m-1)
\geq \left( \frac{c}{s\rho^{1-s}}\right)^m
\left(\frac{n}{m}\right)^{(1-s)m+1}.
\end{split}
\end{equation*}

\vspace{0.2cm}
\vspace{0.3cm}
\begin{theorem}\label{s-energy-RS-LB}
$\!\!\! .\,\,$
There exist reversible averaging systems with initial diameter equal to 1
whose $s$-energy is at least
$\bigl( c/s\rho^{1-s} \bigr)^m
(n/m)^{(1-s)m+1}$,
for constant $c>0$.
The number of vertices is $n$ and the number of connected
components is bounded by $m$; furthermore, all positive entries
in the stochastic matrices are at least $\rho$, with $0<\rho<1/4$.
\end{theorem}
\smallskip
\vspace{0.3cm}

\subsection{Expanding averaging systems: Theorem~\ref{s-Energ-Bounds}~(c)}\label{thmc}

A $d$-regular expander with 
Cheeger constant $h$ is a graph of degree $d$
such that, for any set $X$ of at most half the vertices,
we have $|\partial X|\geq h|X|$, where $\partial X$
is the set of edges with exactly one vertex in $X$.
We say that $G= (V,E)$ is a $d$-regular {\em $m$-expander}
if it has at most $m$ connected components.
Recall that an expanding averaging system {\sc exp}
is a {\sc rev} consisting of $d$-regular $m$-expanders.
Each nonzero entry in $M_t$ is equal to~1 and $q=  d\mathbf{1}$.
We begin the proof of Theorem~\ref{s-Energ-Bounds}~(c)
with a lower bound on the Dirichlet form that exploits the expansion 
of a $d$-regular expander with Cheeger constant $h$.
This is known as {\em Cheeger's inequality}.
We include the proof below for completeness.

\begin{lemma}\label{dirichlet-Expanders}
$\!\!\! .\,\,$
If $G= (V,E)$ is connected, then 
$\sum_{(i,j)\in E}\, (x_i- x_j)^2  \geq b(h/d)^2 \, \|x-\hat x\|_q^2
\geq b(h\Delta)^2/2d$, for constant $b>0$, where 
$\Delta$ is the diameter of the agent positions $x_1,\ldots, x_n$.
\end{lemma}
\proof
All of the ideas in this proof come from~\cite{alon86, sinclair93}.
The inequality is invariant under shifting and scaling, so we may assume
that $\hat x=\mathbf{0}$ and $\|x\|_2= 1$.
Relabel the coordinates of $x$ so they appear in nonincreasing order,
and define $y\in \mathbb{R}^n$ such that $y_i=\max\{x_i,0\}$.
Let $\alpha= \mathrm{argmax}_k(y_k > 0)$ and
$\beta= \min\{\alpha, \lfloor n/2\rfloor\}$. 
By switching $x$ into $-x$ 
if necessary,\footnote{Intuitively, by changing all signs if necessary,
we force the minority sign among 
the coordinates of $x$ to be positive unless their
contribution to the norm of $x$ is too small.}
we can  always assume that
$\|y\|_2^2> c \! := 1/6$ if $\alpha=\beta$, and
$\|y\|_2^2\geq 1-c$ if $\alpha>\beta$.
By Cauchy-Schwarz, $(y_i+y_j)^2\leq 2(y_i^2+y_j^2)$; hence, 

\begin{equation}\label{yiyjUB1}
\begin{split}
\sum_{(i,j)\in E} \, \bigl|\, y_i^2-y_j^2 \,\bigr| 
&= 
\sum_{(i,j)\in E} \, (y_i+y_j) |y_i-y_j|
\leq \sqrt{ \sum_{(i,j)\in E} (y_i + y_j)^2 \sum_{(i,j)\in E} (y_i - y_j)^2} \\
& \leq    \sqrt{ \sum_i 2dy_i^2 \sum_{(i,j)\in E} (y_i- y_j)^2  }
\leq   \sqrt{ \sum_{(i,j)\in E} 2d(y_i- y_j)^2}
\leq   \sqrt{ \sum_{(i,j)\in E} 2d(x_i- x_j)^2}.
\end{split}
\end{equation}
By the expansion property of $G$,
summation by parts yields
\begin{equation}\label{yiyjLB1}
\sum_{(i,j)\in E}\, \bigl|\, y_i^2-y_j^2 \,\bigr| \geq 
\sum_{k=1}^{\lfloor n/2\rfloor} hk \bigl( y_k^2- y_{k+1}^2 \bigr) 
+
\sum_{k=\lfloor n/2\rfloor +1}^{n-1} h(n-k) \bigl( y_{k+1}^2 - y_k^2 \bigr)
= h\, \bigl( \|y\|_2^2- n y_{\lfloor n/2\rfloor +1}^2\bigr).
\end{equation}
Suppose that $\alpha>\beta = \lfloor n/2\rfloor$.
It follows from $\sum_{i=1}^n x_i=0$ that
$\sum_{i=\alpha+1}^n |x_i|= \|y\|_1\geq (\beta+1)y_{\beta+1}$.
By Cauchy-Schwarz, this yields
$n y_{\beta+1}^2/4\leq \sum_{i=\alpha+1}^n x_i^2= 1-\|y\|_2^2$,
and by~(\ref{yiyjLB1}), 
$\sum_{(i,j)\in E}\, \bigl|\, y_i^2-y_j^2 \,\bigr| \geq 
(1-5c)h= c h$.
If $\alpha=\beta$, then $y_{\lfloor n/2\rfloor +1}=0$ and
$\sum_{(i,j)\in E}\, \bigl|\, y_i^2-y_j^2 \,\bigr| \geq h  \|y\|_2^2> ch$.
Applying~(\ref{yiyjUB1}) shows that
$
\sum_{(i,j)\in E} (x_i- x_j)^2 \geq (c h)^2/2d$, which gives us
the first inequality of the lemma. The second one
follows from the fact that the interval $[u,v]$ enclosing
the vertex positions contains~$0$.
By Cauchy-Schwarz, $1= \|x\|_2^2\geq u^2+v^2\geq \frac{1}{2}(v+|u|)^2= \Delta^2/2$,
and the proof is complete.
\hfill $\Box$
\proofend

Let $G_{\leq t}$ be the graph obtained by adding all the edges from
$G_0,\ldots, G_t$.
Let $m_t\leq m$ be the number of connected components in $G_t$,
and let $\Delta_{t,i}$ denote the 
diameter of the $i$-th component of $G_t$ (labeled in any order).
Let $t_1,\ldots, t_c$ be the times $t>0$ at which
the addition of $G_{t}$ reduces the number of components in $G_{\leq t-1}$.
If no such times exist, set $c=1$ and $t_c = 0$.

\begin{lemma}\label{CC-cover-length}
$\!\!\! .\,\,$
If $G_{\leq t_c}$ is connected, then
$\sum_{t\leq t_c} \sum_{i=1}^{m_t} \Delta_{t,i}^2 \geq \, 
\frac{1}{2dmn} \|x-\hat x\, \|_q^2\,$.
\end{lemma}
\proof
At any time $t_k>0$,
the drop $d_k$ in the number of components can be achieved by $d_k$ (or fewer)
components in $G_{t_k}$.
We collect the intervals spanned by these components into a set $F$,
to which we add the intervals for the components of $G_0$; thus $|F|<2m$.
A simple convexity argument (omitted) shows that the
union of the intervals in $F$ coincides with 
the interval $[a,b]$ enclosing the $n$ vertices at time~$1$; so
the lengths $l_1,\ldots, l_{|F|}$ of the intervals in $F$ 
sum up to at least $b-a$.
By Cauchy-Schwarz,
$\sum_{t\leq t_c} \sum_{i=1}^{m_t} \Delta_{t,i}^2 
\geq \sum_{i=1}^{|F|} l_i^2 \geq (b-a)^2/(2m-1)$.
The lemma follows from~$\|x-\hat x\, \|_q^2\leq \|q\|_1(b-a)^2\leq  dn (b-a)^2$.
\hfill $\Box$
\proofend
\smallskip

By~(\ref{DirichletDef}) and Lemma~\ref{dirich},
for any $x= x(0)\in \mathbb{R}^n$,
$\|P_t\, x\|_q^2  \leq \|x\|_q^2 - \frac{D_t}{2}$, so
$\|x\|_q^2 -  \|x(t_c+1)\|_q^2
\geq \frac{1}{2} \sum_{t\leq t_c} D_t
\geq \frac{1}{d} 
\sum_{t\leq t_c} \sum_{(i,j)\in E_t}\, \bigl(x_i(t)- x_j(t) \bigr)^2$.
Assuming that $G_{\leq t_c}$ is connected, 
Lemmas~\ref{dirichlet-Expanders} and ~\ref{CC-cover-length}
imply that 
\begin{equation}\label{Arecur}
\|x\|_q^2 -  \|x(t_c+1)\|_q^2
\geq \frac{bh^2}{2d^2} \sum_{t\leq t_c}\sum_{i=1}^{m_t} \Delta_{t,i}^2
\geq \frac{bh^2}{4d^3mn} \|x-\hat x\, \|_q^2
 \, .
\end{equation}

\noindent
Let $V(n,m)$ be the maximum $s$-energy of an expanding averaging system {\sc exp}
with at most $n$ vertices and $m$ connected components at any time,
subject to the initial condition $\|x - \hat x\, \|_q^2 = 1$ and,
without loss of generality, $\hat x= \mathbf{0}$.
By~(\ref{Arecur}),
$\|x(t)\|_q^2$ shrinks by at least a factor of $\alpha:= 1- bh^2/(4d^3mn)$
by time $t_c+1$.
By scaling, we see that the $s$-energy expanded after $t_c$ 
is at most $\alpha^{s/2} V(n,m)$.
While $t<t_c$ (or if $G_{\leq t_c}$ is not connected),
the system can be decoupled into two {\sc exp} systems
with fewer than $m$ components.
Since $\|x\|_q=1$, the diameter of the system is at most 
$2\max_i |x_i|\leq 2/\sqrt{d}$; therefore
$$
V(n,m) \leq \alpha^{s/2} V(n,m) + 2 V(n,m-1) + m\bigl(2/\sqrt{d}\, \bigr)^s.
$$
It follows that 
\begin{equation}\label{V-recur}
V(n,m) \leq \frac{2}{1-\alpha^{s/2}}\, \Bigl( V(n,m-1) + m \, \Bigr).
\end{equation}
If $m=1$ then $t_c=0$, so we can bypass Lemma~\ref{CC-cover-length}
and its reliance on the diameter. Instead, we use the connectedness
of the graphs to derive from Lemma~\ref{dirichlet-Expanders}:
\begin{equation*}
\|x\|_q^2 -  \|x(1)\|_q^2
\geq \frac{1}{d} \sum_{(i,j)\in E_0}\, (x_i- x_j)^2
\geq \frac{bh^2}{d^3} \, \|x-\hat x\|_q^2
 \, .
\end{equation*}
Setting $\alpha= 1 - bh^2/d^3$ shows that $V(n,m)= O(d^3/s h^2)$ 
for $m=1$. By~(\ref{V-recur}), we find that
$V(n,m) \leq (c/s)^m$, for $c= O(d^3mn/h^2)$,
which completes the proof
of Theorem~\ref{s-Energ-Bounds}~(c).
\hfill $\Box$
\proofend

\subsection{Random averaging systems: Theorem~\ref{s-Energ-Bounds}~(d)}\label{thmd}

It is assumed here that each graph $G_t$ is picked independently, uniformly
from the set of simple $(d-1)$-regular graphs 
with $n$ vertices, with fixed $d>3$~\cite{worlmald99}.
(We use $d-1$ because $d$ must account for the self-loops.)
The system is a special case of a {\sc rev}, so we use the same notation.
The stochastic matrix $P_t$ for $G_t$ is $\frac{1}{d} M_t$,
where $M_t$ is a random symmetric 0/1 matrix
with a positive diagonal and all row sums equal to $d$;
we have $q= d\mathbf{1}$.
Recall that the $s$-energy $\mathcal{E}_s^{\text{\sc ran}}$
is defined over systems with unit variance $\|x-  \hat x\, \|_q^2 =1$, where
$\hat x= \frac{1}{d} \langle x, \mathbf{1}\rangle$.
Let $x_i, y_i$ be the positions of agent $i$ at step $t$ and $t+1$ respectively.
As usual, we may place the center of gravity $\hat x$ at
the origin at time~0, where it will remain forever; that is, $\sum_{i=1}^n x_i= 0$.

\begin{lemma}\label{lemma:step} 
$\!\!\! .\,\,$
$\mathbb{E} \sum_{i=1}^n y_i^2  = (1-b) \sum_{i=1}^n x_i^2$,
where $b= \frac{(d-1)n}{d(n-1)}\,$.
\end{lemma}
\proof
Write $\delta_{ij} = x_i - x_j$
and $M_t= (m_{ij})$.  With all sums extending from $1$ to $n$, we have
\begin{equation}\label{DiffXY}
\begin{split}
\sum_{i} x_i^2 - \sum_{i} y_i^2 
&= \sum_{i} x_i^2 - \sum_{i} \left(x_i - \frac{1}{d}\sum_j m_{ij}\delta_{ij}\right)^2 \\
&= \frac{2}{d}\sum_{i,j}m_{ij}x_i\delta_{ij} - \frac{1}{d^2} \sum_{i,j,k} m_{ij}m_{ik}\delta_{ij}\delta_{ik}\\
&= \frac{1}{d}\sum_{i,j}m_{ij}\delta_{ij}^2 - \frac{1}{2d^2} \sum_{i,j,k} m_{ij}m_{ik}
(\delta_{ij}^2 + \delta_{ik}^2 - \delta_{jk}^2)\\
&= \frac{1}{2d^2}\sum_{i,j,k: i\neq j} m_{ik}m_{jk}\delta_{ij}^2 \, ,
\end{split}
\end{equation}
with the last equality following from $\sum_{k = 1}^n m_{ik} = d$
and $\delta_{ii}= 0$.
By symmetry, $\text{Pr}\, [m_{ij}=1] = (d-1)/(n-1)$ and
$\text{Pr}\, [m_{ij}m_{ik}=1] = \binom{d-1}{2}/  \binom{n-1}{2}$,
for any pairwise distinct $i, j, k$. 
For any $i\neq j$, we have
$\sum_k m_{ik}m_{jk}= 2m_{ij} +  \sum_{k: k\neq i,k\neq j} m_{ik}m_{jk}$;
hence 
$$
\sum_{k}  \mathbb{E}\, [ m_{ik}m_{jk} ] 
= 
2\, \mathbb{E}\, [ m_{ij} ] +
\sum_{k: k\neq i,k\neq j} \mathbb{E}\, [ m_{ik}m_{jk} ]
=   \frac{d(d-1)}{n-1} .
$$
Since $\sum_i x_i=0$, we have
$\sum_{i,j: i\neq j}  \delta_{ij}^2 = 2n \sum_i x_i^2$.
By~(\ref{DiffXY}), it follows that
\begin{equation*}
\mathbb{E} \sum_{i=1}^n y_i^2 
=  \sum_{i} x_i^2 - \frac{1}{2d^2}\sum_{i,j,k: i\neq j} \delta_{ij}^2 \,  \mathbb{E}\, [ m_{ik}m_{jk} ] \\
=  \sum_{i} x_i^2 -
    \frac{d-1}{2d(n-1)}\sum_{i,j: i\neq j}  \delta_{ij}^2 \, .
\end{equation*}
\hfill $\Box$
\proofend

Markov's inequality tells us that
$\sum_i y_i^2 \geq (1-b/3)\sum_i x_i^2$ holds
with probability at most
$\mathbb{E}\, \bigl[\sum_i y_i^2 \bigr] \bigm/ \bigl[(1-b/3)\sum_i x_i^2 \bigr]
\leq 1-b/2$.
Since $\|x\|_q=1$, the diameter of the system is at most $2/\sqrt{d}$;
by the usual scaling law, it follows that
\begin{equation*}
\mathbb{E}\, \mathcal{E}_s^{\, \text{\sc ran}}
\leq  2^s\, \mathbb{E} K + \frac{b}{2}\bigl(1 - b/3\bigr)^{s/2}\,
\mathbb{E}\, \mathcal{E}_s^{\, \text{\sc ran}} +
\bigl(1 - b/2\bigr)\, \mathbb{E}\, \mathcal{E}_s^{\, \text{\sc ran}},
\end{equation*}
where $K$ is the number of connected components in $G_0$.
It is known~\cite{worlmald99} 
that, for $d>3$, the probability that the graph is
not connected is $O(n^{3-d})$; hence $\mathbb{E} K = O(1)$.
Since $b\geq 1/2$, we conclude that
$\mathbb{E}\, \mathcal{E}_s^{\, \text{\sc ran}}
= O\bigl(1/\bigl(1- (5/6)^{s/2} \bigr)\bigr)=  O(1/s)$;
hence Theorem~\ref{s-Energ-Bounds}~(d).
\hfill $\Box$
\proofend

\section{Applications}\label{app1}

We use Theorem~\ref{ConvergeUB} to bound 
the convergence rates of systems for bird flocking, opinion dynamics, 
and distributed motion coordination. We also discuss
the concept of the ``Overton Window'' and provide a theoretical validation for it.

\begin{itemize}
\item
In~\S\ref{birdFlock}, we establish sufficient conditions for the quick
relaxation to kinetic equilibrium in the classic Vicsek-Cucker-Smale model of bird flocking~\cite{chazFlockPaperI, CuckerSmale1, vicsekCBCS95}.
The convergence time is polynomial in the number of birds
as long as the number of flocks remains bounded.
This new result relies on Theorem~\ref{ConvergeUB} as well as 
novel insights into the convex geometry of flocking.
\item
In~\S\ref{PatternFormation}, we investigate a distributed motion coordination
algorithm introduced by~Sugihara and Suzuki~\cite{sugihara-suzuki-1990, bulloBk}.
The idea is to use a swarm of robots to produce a preset pattern, in this case
a polygon. We prove a polynomial bound on the relaxation time of this process.
We enhance the model by allowing faulty communication and proving
that the end result is robust under stochastic errors. We also generalize
the geometry to 3D and arbitrary communication graphs.
\item
In~\S\ref{overton}, 
we apply Theorem~\ref{ConvergeUB} to 
opinion formation in social networks.  We extend the model
to include directed edges so as to capture both evolving and fixed sources
of information. We show that, while all opinions 
might keep changing forever, they will inevitably land
in the convex hull of the fixed sources.  Furthermore, we 
bound the time at which this must happen. 
Our result provides a quantitative validation 
of the {\em Overton window}
as an attracting manifold of ``viable'' 
opinions~\cite{bevensee, bobric, dyjack20, fedyanin, morgan20}.
\end{itemize}

\subsection{Bird Flocking}\label{birdFlock}

Introduced by Reynolds~\cite{reynolds87} in 1987, three heuristic rules have been used widely to
produce spectacular bird flocking animations. 
The three flocking rules are (1) \textit{separation}: avoid collision
(2) \textit{cohesion}: stay grouped together,
and (3) \textit{alignment}: align headings.
Several models are constructed based on these rules to understand flocking dynamics.

We study a variant of the classic Vicsek-Cucker-Smale 
model~\cite{CuckerSmale1, vicsekCBCS95}, a group of $n$ birds are flying in 
the air while interacting via a time-varying 
network~\cite{blondelHOT05, chazFlockPaperI, HendrickxB, jadbabaieLM03}.
The vertices of the network correspond to the
$n$ birds and any two birds are joined by an edge if their distance is at most some 
fixed $r\leq 1$.
The flocking network $G_t$ is thus undirected.
Its connected components define the {\em flocks}.
Each bird $i$ has a position $x_i(t)$ and a velocity $v_i(t)$,
both of them vectors in $\mathbb{R}^3$. Given the state of the system
at time $t=0$, we have the recurrence: for any $t\geq 0$, 

\begin{equation}\label{modelVCS}
\begin{cases}
\, x_i(t+1)= x_i(t)+ v_i(t+1); \\
\, v_i(t+1)= v_i(t) + a_i\sum_{j\in \mathcal{N}_i(t)}\bigl( v_j(t) - v_i(t)\bigr),
\end{cases}
\end{equation}
where $\mathcal{N}_i(t)$ is the set of vertices $j\neq i$ adjacent to $i$ at time $t$.
At each step, a bird adjusts its velocity
by taking a weighted average with its neighbors.
The weights $a_i$ 
indicate the amount of influence birds
exercise on their neighbors.  To avoid negative weights, we
require that $0<a_i\leq  1/(|\mathcal{N}_i(t)|+1)$. We write $\rho:= \min_i a_i \in (0,1/2]$.

Intuitively, by repeating the recurrence, each bird should eventually converge to 
a fixed speed and direction.
This is supported by computer simulations and several convergence 
results~\cite{HendrickxB, Moreau2005, olfati06}.
As was shown in~\cite{chazFlockPaperI}, however, the model above might be periodic
and never stabilize.  To remedy this, we stipulate that, for
two birds to be newly joined by an edge,  
their velocities must differ by at least a minimum amount:
Formally, we require that, at any time $t$, 
$(i,j)\in G_t\setminus G_{t-1}$ if $\|x_i(t)-x_j(t)\|\leq r$ 
and $\|v_i(t)- v_j(t)\| > \eps_o$, for small fixed positive~$\eps_o$.
By space and time scale invariance, we may assume\footnote{These
bounds are nonrestrictive and the choice of $\frac{1}{2} \sqrt{\rho / n}$ is made only
to simplify some calculations.}
that 
$\|x_i(0)\| \leq 1$ and $\|v_i(0)\| \leq \frac{1}{2}\sqrt{\rho / n}$, for all birds~$i$.
We state our main result. Its 
main novelty is that the convergence time is polynomial in the number of birds,
as long as the number of flocks is bounded by a constant.
\vspace{0.4cm}

\vspace{0.3cm}
\begin{theorem}\label{ConvergeTime}
$\!\!\! .\,\,$
A group of $n$ birds forming a maximum of $m\leq n$ flocks
relax to within $\eps$ of a fixed velocity vector in 
time $O(n^2/\rho)\log (1/\eps) +t_o$, where
$t_o$ is on the order of
$ m n^{2(m+2)}(1/\rho)^{m+1}\log \frac{n}{\rho}\,$.
\end{theorem}
\smallskip

\vspace{0.4cm}

\vspace{0cm}
\begin{figure}[htb]
\begin{center}
\hspace{0cm}
\includegraphics[width=6cm]{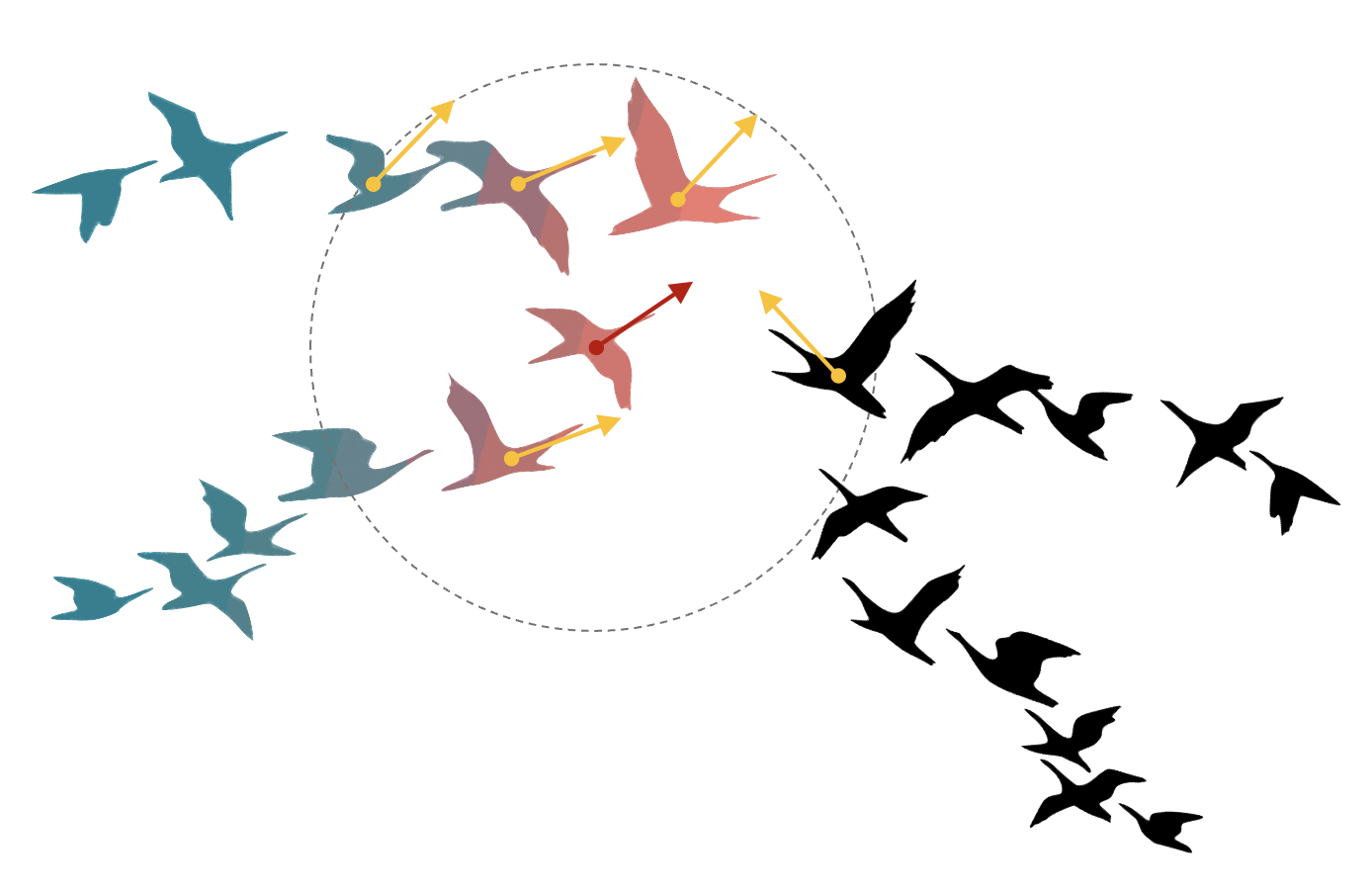}
\end{center}
\vspace{-0.5cm}
\caption{\small A bird is influenced by its neighbors within distance $r$. \label{fig-flock}}
\end{figure}
\vspace{0.5cm}

We rewrite the map of the velocity dynamics~(\ref{modelVCS}) in matrix form,
$v(t+1)= P_t v(t)$, for $t\geq 0$, where $v(t)$ is an $n$-by-3 matrix with each row indicating
a velocity vector. We have
$P_t= Q^{-1}M_t$, where
$Q= \hbox{diag}(q)$,
$q_i= 1/a_i$, and

\begin{equation*}
(M_t)_{ij} =
\begin{cases}
\,\, \text{$1/a_i- |\mathcal{N}_i(t)|$ \ if $i=j$}  \\
\,\, \text{1 \ if $(i,j) \in G_t$ and $i\neq j$} \\
\,\, \text{0 \ else.}
\end{cases}
\end{equation*}
Note that $\bar q:= \|q\|_1^{-1} q$ is the joint stationary distribution
and $q= M_t\mathbf{1}\leq  \mathbf{1}/\rho$, where
$\rho:= \min_i a_i \in (0,1/2]$. 
Each nonzero entry in row $i$ of $P_t$, being either
$a_i$ or $1-a_i |\mathcal{N}_i(t)|\geq a_i$, is at least $\rho$.
This shows that each one of the three coordinates provides
its own reversible averaging system $\mathcal{S}_j$ ($j=1,2,3$).
The only difference between the systems is their initial states.
The $s$-energy of any such system is equal to
$\sum_t E_{s,t}$, where $E_{s,t}= \sum_i l_i(t)^s$
and $l_i(t)$ is the length of the $i$-th block at time~$t$.
(Recall that a block is an interval formed by the embedded edges.)
Let $m$ be the maximum number of flocks and
$N_{m,\alpha}$ the number of times~$t$ at which
some block length $l_i(t)$ from at least one of 
$\mathcal{S}_j$ ($j=1, 2, 3$) is at least $\alpha$.
For $0<\alpha<1$, we have 
$N_{m,\alpha}\leq \inf_{s\in (0,1]} 3 \alpha^{-s} \mathcal{E}_{m,s}^{\text{\sc rev}}$.
Our assumption that $\|v_i(0)\| \leq \frac{1}{2}\sqrt{\rho / n}$ for all birds $i$
implies that each of the three systems has variance at most~1.
Setting $s=c/\log (1/\alpha)$ and assuming that $\alpha$ is small enough so
that $0<s\leq 1$, it follows from
Theorem~\ref{s-Energ-Bounds}~(b) that, for (another) constant $c>0$,

\begin{equation}\label{N-alpha}
N_{m,\alpha} \leq c \Big(\frac{n^2}{\rho}  \log \frac{1}{\alpha} \Big)^m .
\end{equation}

\subsubsection{Single-flock dynamics}\label{SF-dyn}

Between two consecutive switches (ie, edge changes), the flocking networks consists
of fixed non-interacting flocks. We can analyze them separately.
Without loss of generality, assume that $G_t$ is a connected, time-invariant graph.
We focus on system $\mathcal{S}_1$ for convenience.
It consists of a single block at each timestep, so the $s$-energy is
of the form $\sum_t  \Delta_t^s$, where $\Delta_t$ is the diameter of
the system at time $t\geq 0$. 
The diameter can never grow; therefore $\Delta_t \leq \alpha$ for any 
$t\geq N_{1,\alpha}$. By~(\ref{N-alpha}),
it follows that $\Delta_t\leq e^{-a \rho t/n^2}$,
for any $t\geq n^2/a\rho$ and a small enough constant $a>0$.
Recall that $x(t)$ and $v(t)$ are $n$-by-3 matrices; 
denote their first column by $y(t)$ and $w(t)$, respectively.
Write $y(0)= y$, $w(0)=w$, and $P_t=P$.
The vector $w(t) = P^t w$ tends to $(\bar q^T w) \mathbf{1}$.
Since its coordinates lie in an interval of width $\Delta_t$, it follows that
$w(t)= (\bar q^T w) \mathbf{1} + \zeta(t)$, where
$\| \zeta(t) \|_\infty \leq \Delta_t\leq  e^{- a\rho t/n^2}$.
Thus, for some $\gamma , \eta_t \in \mathbb{R}^n$ and any $t\geq  t_a$, with
\begin{equation}\label{ta}
t_a:= (n^2/a\rho)^2,
\end{equation}
we have
$$
y(t) = y + \sum_{k=1}^{t} w(k)= y +  
t (\bar q^T w) \mathbf{1} + \sum_{k=1}^t \zeta(k)
= \beta t + \gamma + \eta_t,
$$
where
$\beta= ( \bar q^T w ) \mathbf{1}$,
$\gamma= y+ \sum_{k=1}^\infty \zeta(k)$,
and $\| \eta_t \|_\infty\leq e^{-b\rho t/n^2}$, for constant $b>0$.
The same holds true for the other two coordinates, so
the birds in the flock fly parallel to a straight line with a deviation
from their asymptotic line vanishing exponentially fast.
If so desired, it is straightforward to lock the flocks by stipulating that 
no two birds can lose an edge between them unless their difference in velocity
exceeds a small threshold~$\theta$; because of the 
exponential convergence rate, choosing $\theta$ small enough
ensures that two birds $i$ and $j$ adjacent in a flock may exceed
distance~$r$ by only a tiny amount.

\subsubsection{Flock fusion}

To bound the relaxation time, we begin with an
intriguing geometric fact:
Far enough into the future, two birds can only come close
to each other if their velocities are nearly identical.
In other words, encounters at large angles of attack 
cannot occur over a long time horizon.
We begin with a technical lemma: A stationary
observer positioned at the initial location of a bird
sees that bird move less and less over time;
this is because the bird flies increasingly in the direction of the line of sight.

\begin{lemma}\label{line-of-sight}
$\!\!\! .\,\,$
There is a constant $c=c(\eps_o)$ such that, for any $i$ and $t>1$,
$$
\Big\| \, v_i(t) - \frac{1}{t} \bigl(x_i(t) - x_i(0)\bigr) \, \Big\|
\leq c n^{2(m+2)}
\left( \frac{1}{\rho}\right)^{m+1}
\left( \frac{\log t}{t} \right) \, .
$$

\end{lemma}
\proof
For notational convenience, we set $i=1$ and
we denote by $y_j(t)$ (resp. $w_j(t)$)
the first coordinate of $x_j(t)$ (resp. $v_j(t)$).
The line-of-sight direction of bird 1 is given by $\frac{1}{t}\bigl( x_1(t)-x_1(0) \bigr)$.
Along the first coordinate axis, this gives
\begin{equation}\label{u=sum}
u:= \frac{1}{t}\bigl( y_1(t) - y _1(0) \bigr)
 =  \frac{1}{t} \sum _{k=1}^{t} w_1(k).
\end{equation}
Consider the difference $\delta:= u - w_1(t)$. We can
define the corresponding quantity for each of the other two directions
and assume that $\delta$ has the largest absolute value among the three of them.
By symmetry, we can also assume that $\delta\geq 0$; therefore
\begin{equation}\label{u-diff}
\bigl\| v_1(t) - \hbox{$\frac{1}{t}$} \bigl(x_1(t)- x_1(0) \bigr) \bigr\|
\leq \sqrt{3}\, \delta. 
\end{equation}

The proof of the lemma rests on showing that, if $\delta$
is too large, some bird $l$ must be at a distance greater than 2
from bird $1$ at time 0, which has been ruled out.
To identify the far-away bird~$l$,
we start with $l=1$ at time $t$, and we trace
the evolution of its flock backwards in time, always trying to move
away from bird~1, if necessary by switching bird $l$ with a neighbor.
This is possible because of two properties, at least one of which holds
at any time $k$:  (i) bird~$l$ flies nearly straight in the time interval $[k, k+1]$;
or (ii) bird~$l$ is adjacent to a bird~$l'$ whose velocity points in a favorable
direction. In the latter case, we switch focus from $l$ to $l'$.

The $s$-energy plays the key role in putting numbers behind these 
properties. For this reason, we define $\mu_l(k)$ as the length of
the block of $\mathcal{S}_1$ containing $w_l(k)$ with respect to 
the flocking network $G_k$. 
Note that $\mu_l(k)$ is the length of an 
interval that contains the numbers $w_j(k)$ for all the birds $j$
in the flock of bird $l$ at time~$k$.
We define the sequence of velocities $\bar{w}(k) = w_l(k)$,
for $k= t, t-1,\ldots, 1$ and $l= l(k)$.
Fix some small $\alpha$ ($0<\alpha \leq \eps_o$).

\vspace{0.1cm}

{\small
\par\medskip
\renewcommand{\sboxsep}{0.5cm}
\renewcommand{\sdim}{0.8\fboxsep}
\begin{center}
\shabox{\parbox{10cm}{
\begin{itemize}
\item[\text{[1]}]
\hspace{0.2cm}
$\bar{w}(t) \leftarrow w_1(t)$ \ \ and \ \ $l \leftarrow 1$
\item[\text{[2]}]
\hspace{0.2cm}
{\bf for }\  $k=t-1, \ldots, 1$
\item[\text{[3]}]
\hspace{0.8cm}
{\bf if } \ $\mu_l(k) > \alpha$
\ {\bf then } \ $l\leftarrow \text{argmin} \left\{\, w_j(k) \, |\,  j\in N_l(k)\,  \right\}$
\item[\text{[4]}]
\hspace{0.8cm}
$\bar{w}(k) \leftarrow w_l(k)$
\end{itemize}
}}
\end{center}
\par
}
\vspace{0.4cm}
Perhaps the best way to understand the algorithm is first to imagine
that the conditional in step [3] never holds:
In that case, $l=1$ throughout and we are simply
tracing the backward evolution of bird 1.  Step [3] aims to catch
the instances where the reverse trajectory inches excessively toward
the initial position of bird 1.
When that happens, $|w_l(k+1)- w_l(k)|$ is large, hence so
is $\mu_l(k)$, and step [3] kicks in.
We exploit the fact that $w_l(k+1)$ is a convex
combination of $\bigl\{w_j(k) \, |\,  j\in N_l(k)\bigr\}$ to update
the current bird $l$ to a ``better'' one.
Using summation by parts, we find that
\begin{equation}\label{sum-parts}
 \sum_{k=1}^{t} \bar{w}(k) =
 t \bar w(t) -
 \sum_{k=1}^{t-1} k\bigl( \bar{w}(k+1) - \bar{w}(k)\bigr).
 \end{equation}
 
Let $R$ be the set of times $k$ that pass the test in step [3]
and $S$ the set of switches (ie, network changes).
An edge creation entails a block of length $\eps_o/\sqrt{3}$ or more
in at least one of $\mathcal{S}_j$ ($j=1,2,3$).
The steps witnessing edge deletions outnumber those 
seeing edge creations by at most a factor of $\binom{n}{2}$.
Pick any two consecutive switches and let $I$ be the time interval between them.
Each flock remains invariant during $I$; thus $|R\cap I|\leq N_{1,\alpha}$;
hence 
\begin{equation}\label{SR}
|S|\leq  n^2 N_{m, \eps_o/\sqrt{3}} 
\hspace{0.7cm} 
\text{and}
\hspace{0.7cm}
|R|\leq  (N_{1,\alpha} +1) |S|.
\end{equation}
\noindent
Because of the single-flock invariance, the diameter of $\mathcal{S}_1$ during $I$ 
can never increase; therefore $J= I\setminus R$ consists of a single time interval.
If $k\in J$, then  
$| \bar{w}(k+1) - \bar{w}(k) |
= |w_l(k+1) - w_l(k) | \leq \mu_l(k)\leq \alpha$ and, 
by~Theorem~\ref{s-Energ-Bounds}~(b),
$
\sum_{k\in J}  | \bar{w}(k+1) - \bar{w}(k) |
\leq \sum_{k\in J} E_{1,k}\leq \alpha \, \mathcal{E}_{1,1}^{\text{\sc rev}}
= O(\alpha n^2/\rho)$;
hence

\begin{equation}\label{Ek-alpha}
\sum_{k\in \{1,\ldots, t-1\}\setminus R}  \bigl| \bar{w}(k+1) - \bar{w}(k)  \bigr| =
O(\alpha n^2|S|/\rho).
\end{equation}

\noindent
Let $l'$ be the value of $l$ in the final assignment
$\bar{w}_l(1) \leftarrow w_l(1)$ in step [4].
Since $\bar{w}(k+1) \geq \bar{w}(k)$ for $k\in R$
and $\bar{w}(t)= w_1(t)$, it follows 
from~(\ref{sum-parts}, \ref{Ek-alpha}) and $r\leq 1$ that
\begin{equation}\label{ylo-y}
\begin{split}
y_{l'}(0) - y_1(0)
&\geq  
\bigl(y_1(t) - y_1(0)\bigr) + \bigl(y_{l'}(0) -y_1(t) \bigr) 
\geq  tu -  \sum_{k=1}^{t} \bar{w}(k) - r |R|  \\
& \geq  t \delta 
+  \sum_{k=1}^{t-1} k\bigl( \bar{w}(k+1) - \bar{w}(k)\bigr)
      - |R|  \geq  t \delta -  O(t \alpha n^2 |S| /\rho)  -  |R|.
\end{split}
\end{equation}

\noindent
We set $\alpha= \eps_o /t$.
Noting that $y_{l'}(0) - y_1(0) \leq 2$,
the lemma follows from~(\ref{N-alpha}, \ref{u-diff}, \ref{SR}) and
$$
\delta\leq d n^{2(m+2)}
\left( \frac{1}{\rho}\right)^{m+1}
\left( \log \frac{1}{\eps_o} \right)^m 
\frac{\log (t/\eps_o)}{t} \, .$$
for constant $d>0$.
\hfill $\Box$
\proofend

\subsubsection{Stabilization}

By Lemma~\ref{line-of-sight}, for a large enough constant $C= C(\eps_o)$,
after time $t>t_0$, where
$$
 t_o:= C m n^{2(m+2)} \left(\frac{1}{\rho}\right)^{m+1} \log \frac{n}{\rho},
$$
no bird's velocity differs from its line-of-sight vector
$\omega_i= \frac{1}{t}\bigl( x_i(t) - x_i(0) \bigr)$ by a vector longer than $\eps_o /3$.
Suppose that birds $i$ and $j$ are within distance $r$ of each other.
By the triangular inequality, 
$ \|\omega_i - \omega_j\| \leq \frac{1}{t}\| x_i(t) - x_j(t)\| + \frac{1}{t} \|x_i(0) - x_j(0)\|\leq (1+r)/t$;
therefore,
$$
\|v_i(t)- v_j(t)\|\leq \|v_i(t)- \omega_i\| + \|\omega_i - \omega_j\| +  \|v_j(t)- \omega_j\| \leq \eps_o.
$$
This implies that each flock is time-invariant past time $t_o$.
By~(\ref{ta}), $t_o\geq t_a$, so
the birds within each flock align their velocities exponentially fast 
from that point on. Theorem~\ref{ConvergeTime} follows.
\hfill $\Box$
\proofend

\subsection{Distributed Motion Coordination}\label{PatternFormation}

In \cite{sugihara-suzuki-1990}, Sugihara and Suzuki introduced 
an interesting model of pattern formation in a swarm of robots. In their model,
the robots can communicate anonymously and adjust their positions
accordingly. Assume that their goal is to align themselves along
a line segment $ab$. Two robots position themselves manually at the endpoints 
of the segment while the others attempt to reach $ab$ by linking with their right/left
neighbors and averaging their positions iteratively.   This setup creates
a polygonal line $u_1=a, u_2,\ldots, u_{n-1}, u_n=b$, where $u_i$ is
the position $(x_i,y_i)$ of robot~$i$. The polygonal line converges to $ab$
in the limit. We use the $s$-energy bounds
to evaluate the convergence time of the robots. We actually prove 
a stronger result by generalizing the model in two ways:
(i) we consider the case of an arbitrary communication 
network of robots in 3D, with a subset of vertices pinned to a fixed plane;
(ii) the network suffers from stochastic edge failures.
Our model trivially reduces to Sugihara and Suzuki's by projection.
Allowing stochastic failures to their motion coordination model is novel. 

Let $G$ be a connected (undirected) graph with $n$ vertices labeled in $[n]$
and no self-loops,
and let  the {\em communication weights} $a_1,\ldots, a_n$ be $n$ positive reals
such that $a_i<1/(d_i+1)$,
where $d_i$ is the degree of vertex $i$. We define
$d= \max d_i$ and $\rho= \min a_i$.
For any $t\geq 0$, we define $G_t$ by deleting each edge of $G$ with
probability $1-p$.
We define a (random) stochastic matrix $P_t$ for $G_t$ as follows:
\begin{quote}
\begin{enumerate}
\item
Initialize $P_t=0$;
\item
If $(i,j)$ is an edge of $G_t$, we set $(P_t)_{ij} = a_i$ and $(P_t)_{ji}=a_j$.
\item
$(P_t)_{ii}= 1- \sum_{j (j\neq i)} (P_t)_{ij}$, for all $i$.
\end{enumerate}
\end{quote}
\vspace{0.1cm}
Note that every positive entry of $P_t$ is at least $\rho$.
We embed $G$ in $\mathbb{R}^3$ and pin a subset $R$ of $r$ vertices
to a fixed plane. We fix the scale by assuming that
the embedding lies in the unit cube $[0,1]^3$. 
Without loss of generality, we choose the plane $X=0$.
To ensure the immobility of the $r$ vertices,
we use {\em symmetrization}~\cite{chazelle-Energ1-2011}.
This entails attaching to $R$ a copy 
of $G$ and initializing the embedding of the two copies as mirror-image
reflections about $X=0$.
This modification increases the number
of vertices to $\nu= 2n-r$.
The sequence $(P_t)_{t\geq 0}$ is defined by picking a random $G_t$ (as defined above)
at each step~{\rm iid}.

The vertices of $R$ are embedded in the plane $X=0$ at time 0, where,
by symmetry, they reside permanently (but may move within the plane).
To assert and measure the attraction of 
the points to the plane, it suffices to focus on the dynamics along
the $X$-axis.  Given $x(0)\in [-1,1]^\nu$, we have
$x(t+1) = P_t x(t)$. This gives us a reversible averaging system,
with $q= (1/a_1,\ldots, 1/a_\nu)$; note that symmetrization may increase
the degrees by at most a factor of 2; and so $a_i$ and hence $\rho$ may
have to be scaled down by up to one half. 
Put $\|y\|_q:= \sqrt{ \langle y, y \rangle_q}$, where
$\langle y, z \rangle_q : = \sum_i q_i y_i z_i$ and
$$D_t= \sum_{i}\max_{j:\, (i,j)\in G_t} \bigl( x_i(t)- x_j(t) \bigr)^2 . $$
Since $G$ is connected, there is a path $\pi$ connecting the leftmost
to the rightmost vertex along the $X$-axis.
By the Cauchy-Schwarz inequality,
\begin{align*}
\mathbb{E} \, D_t
\geq \mathbb{E} \sum_{i=1}^\nu  \max_{j:(i,j)\in G_t}\delta_{ij}^2
&\geq   \sum_{i=1}^\nu \sum_{j:(i,j)\in G} p \delta_{ij}^2/d_i \\
\geq \frac{p}{d\nu} \Bigl(\sum_{(i,j)\in \pi} |\delta_{ij}| \Bigr)^2
&\geq  \frac{\rho p}{d\nu^2} \|x\|_q^2 \, ,
\end{align*}
\noindent
where $\delta_{ij} = \delta_{ij}(t) = x_i(t) - x_j(t)$.
It follows from Lemma~\ref{dirich} (see~\S\ref{RAS} below) that, 
for $c:= \rho p / (2d \nu^2)$,
$$
\mathbb{E} \, \|Px\|_q^2  \leq \|x\|_q^2 - \frac{1}{2}\,  \mathbb{E} \, D_t
\leq ( 1 - c ) \|x\|_q^2.
$$
By Markov's inequality,
$$
\Pr\left[\|Px\|_q^2 \geq \Bigl(1-\frac{c}{3}\Bigr) \|x\|_q^2\right] 
\leq \frac{\mathbb{E}\,  \|Px\|_q^2  }{ (1 - c/3)  \|x\|_q^2}\leq
1-\frac{c}{2} \, .
$$

\begin{figure}[htb]
\begin{center}
\hspace{0cm}
\subfloat[Initially, moving robots are placed randomly in the unit cube.]{
\includegraphics[width=7cm]{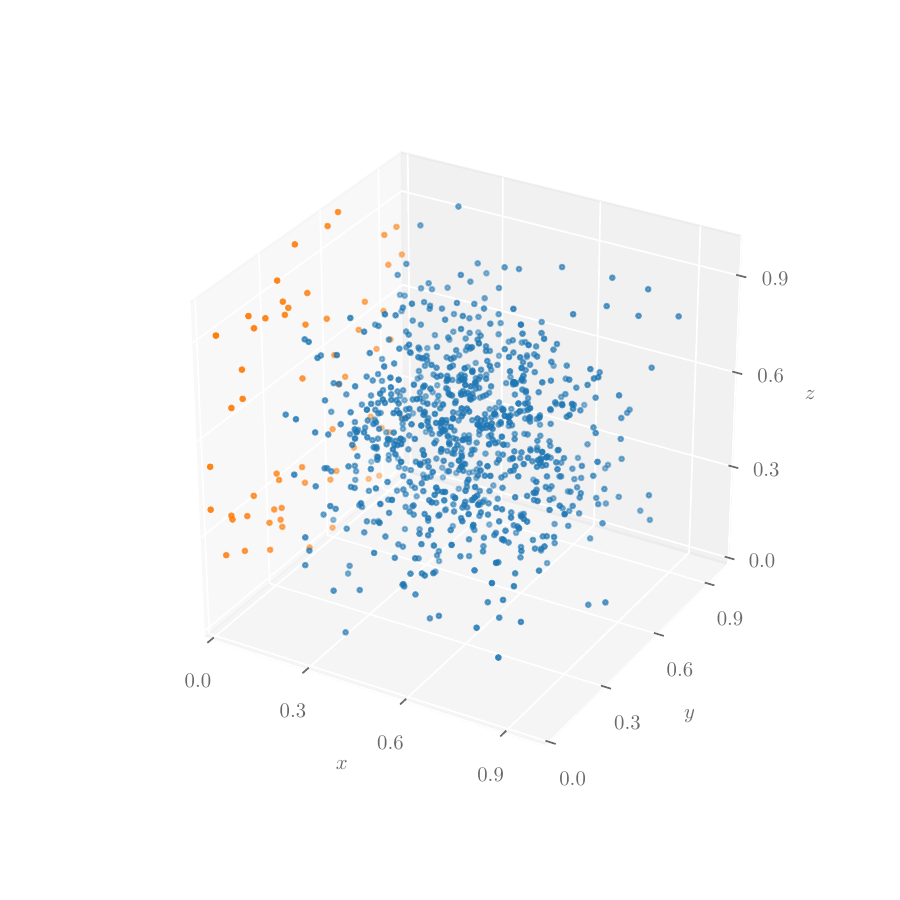}
}
\hspace{0cm}
\subfloat[Positions of moving robots after 200 steps.]{
\includegraphics[width=7cm]{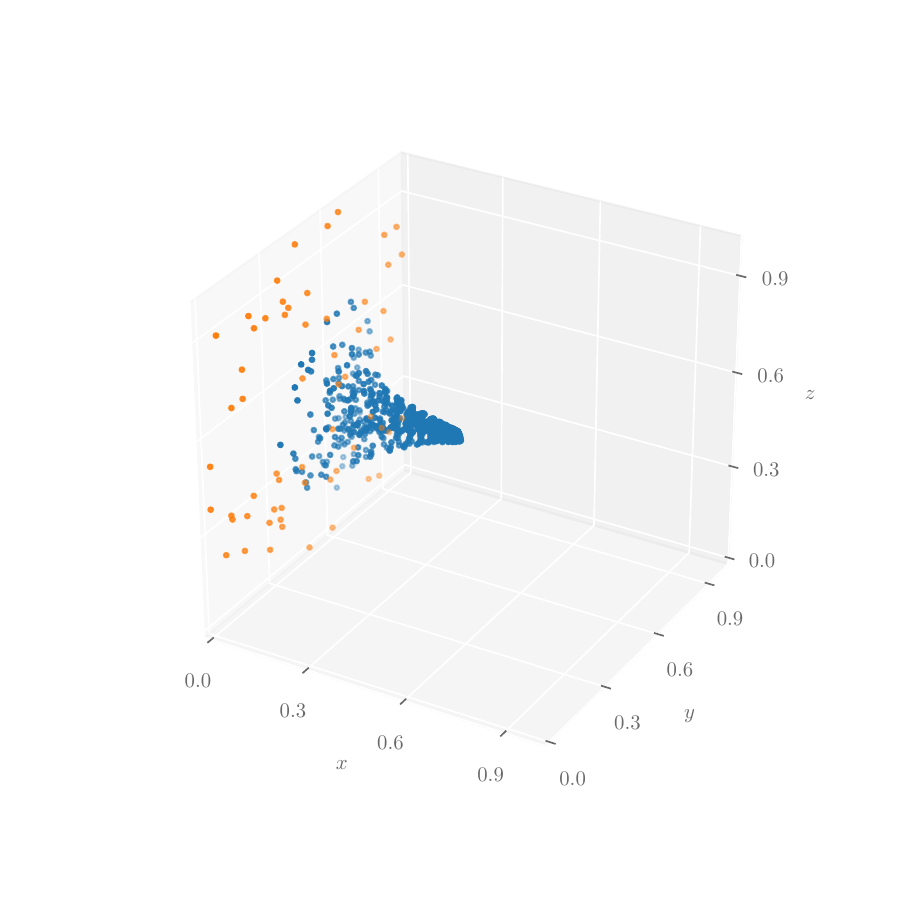}
}
\end{center}
\vspace{-0.5cm}
\caption{\small 
Simulation of a network of 900 robots in 3D with 60 of them (orange dots) pinned to
a fixed plane $X=0$ and the rest (blue dots) are moving. The underlying network $G$ is a 30-by-30 grid graph, 
with edge failure probability equal to 0.3. 
The 60 nodes on two opposite sides of the grid are pinned to $X=0$.
\label{fig-robots}}
\vspace{1cm}
\end{figure}
\vspace{0cm}

\noindent 
Let $l_1(t),\ldots, l_{k_t}(t)$ be the lengths of the blocks
formed by the edges of $G_t$ embedded along 
the $X$-axis.\footnote{Recall that the blocks are the intervals formed
by the union of the embedded edges of $G_t$.}
In a slight variant, we define the $s$-energy $E_s= \sum_{t\geq 0} E_{s,t}$,
where $E_{s,t}= \max_{i=1}^{k_t} l_i(t)^s$.
We denote by $W_s$ the maximum expected 
$s$-energy, where the maximum is taken over all initial positions
such that $\|x\|_q^2 \leq \nu/\rho$.
Since the vertices are embedded in $[-1,1]$ with symmetry about the origin,
this applies to the case at hand.  
Since $q_i= 1/a_i>1$, 
the initial diameter for $W_s$ is at most
$\sqrt{2}\, \|x\|_q \leq \sqrt{2\nu/\rho}$. (Note that we cannot use~2
as the bound for the diameter since $W_s$ only bounds the variance and
not the diameter.)
By scaling invariance, we have the following recurrence relation:
\begin{equation*}
\begin{split}
W_s &\leq \left(\frac{2\nu}{\rho}\right)^{s/2} + 
\frac{c}{2}\left(1 - \frac{c}{3}\right)^{s/2}W_s
+ \left(1 - \frac{c}{2}\right)W_s
\leq  \frac{2^{s+1}\nu^{s/2}}{c\rho^{s/2}\big(1 - (1-c/3)^{s/2}\big)}\\
&= O\left(\frac{\nu^{s/2}}{sc^2\rho^{s/2}}\right)
=  O\left( \frac{d n^2}{\rho p}\right)^2 \frac{(n/\rho)^{s/2}}{s} \, .
\end{split}
\end{equation*}
Let $N_\alpha$ be the number of times $t$ at which some block length $l_i(t)$ 
is at least $\alpha$. 
For $0 < \alpha < 1$, we have $\mathbb{E}\, N_\alpha  \leq \inf_{s \in (0,1]} \alpha^{-s}W_s$. 
Setting $s = 1/\log (n/\rho \alpha^2)$ yields 
$$
\mathbb{E}\,N_\alpha =  
O\left( \frac{d n^2}{\rho p}\right)^2 \log\frac{n}{\rho\alpha}\, .
$$
Let $K_\alpha$ be the number of times $t$ at which there exists an edge $(i,j) \in G_t$ 
whose length $|\delta_{ij}(t)|$ exceeds $\alpha$;
obviously, $K_\alpha \leq N_\alpha$.
Let $T_\alpha$ be the last time at which 
the diameter of the system exceeds $\alpha$.
For each $t \leq T_\alpha$, being a connected graph,
$G$ must include an edge $(i,j)$
whose length $|\delta_{ij}(t)|$ exceeds $\alpha/\nu$. 
That edge belongs to $G_t$ with probability $p$;
therefore $\mathbb{E}\, K_{\alpha/\nu} \geq p \,\mathbb{E}\, T_\alpha$;
hence $\mathbb{E}\, T_\alpha \leq \frac{1}{p} \mathbb{E}\, N_{\alpha/\nu}$.

\vspace{0.2cm}
\begin{theorem}\label{convergence-bound-DMC}
$\!\!\! .\,\,$
The robots align themselves within distance $\epsilon<1$ of a fixed plane  
in expected time $O\bigl(d^2n^4/p^3 \rho^2 \bigr) \log(n/\rho\eps)$,
where $d$ is the maximum degree of the underlying communication network,
$n$ is the number of robots, $1-p$ is the probability of edge failure,
and $\rho$ is the smallest communication weight.
\end{theorem}
\vspace{0.2cm}

\subsection{The Overton Window Attractor}\label{overton}

Following in a long line of opinion dynamics 
models~\cite{castellanoFL2009, deffuant, degroot1974, fj1990, hegselmanK},
we consider a collection of
$n$ agents, each one holding an opinion vector $x_i(t)\in [0,1]^d$ at time $t$;
we denote by $x(t)$ the $n$-by-$d$ matrix whose $i$-th row corresponds to $x_i(t)$.
Given a stochastic matrix $P_t$, the agents update their opinion vectors
at time $t \geq 0$ according to the evolution equation $x(t+1)= P_t \, x(t)$.   
We assume that the last $k$ agents $n-k+1,\ldots, n$ are {\em fixed}
in the sense that $x_i(t)$ remains constant at all times $t \geq 0$.
Algebraically, the square block of $P_t$ corresponding
to the $k$ fixed agents is set to the identity matrix $\mathbb{I}_k$.
The fixed agents can influence the mobile ones, but not the other way around.
The presence of fixed agents (also called ``stubborn,'' ``forceful'' or ``zealots'' in the
literature) has been extensively 
studied~\cite{acemoglu-misinfo, abrahamsson-random, ChazW17,
ghaderi2014,
mobilia2003, mobilia2007, tian-stubborn,
yildiz-binary}. 

In the context of social networks,
the fixed sources may consist of venues with low user influence, such
as news outlets, wiki pages, influencers, TV channels, political campaign sites, 
etc.~\cite{chitra-musco, gaitonde21, JWWSHDZ-latent-online,
lewisGK, NHHTK-filter-bubble, VMCG-evolution, zafaraniAL}.
We know how the mobile agents migrate to the convex hull of the
fixed agents; crucially, we bound the rate of attraction. 
This provides both a quantitative illustration of the famous {\em Overton window}
phenomenon as well as a theoretical explanation for why the window
acts as an attracting manifold~\cite{bevensee, bobric, dyjack20, fedyanin, morgan20}.
Interestingly, the emergence of a global attractor does not imply convergence
(ie, fixed-point attraction). The mobile agents might still fluctuate widely
in perpetuity. The point is that they will always do so within the confines of the
global attractor.

To reflect the stochasticity inherent in the choice of sources visited by a user
on a given day, we adopt a classic ``planted'' model:
Fix a connected $n$-vertex graph $G$ and two parameters
$p\in (0,1]$ and $\rho \in (0, 1/2]$.
At each time $t \geq 0$, $G_t$ is defined by picking every edge of $G$ with
probability at least $p$. (No independence is required
and $n$ self-loops are included.)
We define an $n$-by-$n$ stochastic matrix $P_t$ 
by setting every entry to 0 and updating it as follows:
\begin{enumerate}
    \item For $i > n-k$, $(P_t)_{ii} = 1$.
    \item For $i\leq n-k$, set $(P_t)_{ij} \geq \rho$ for any $j$ such that 
    $(i,j)$ is an edge of $G_t$.
\end{enumerate}
Note that the update is highly nondeterministic. The only two conditions
required are that (i) nonzero entries be at least $\rho$ and (ii) each row sum up to~1.

\vspace{0.5cm}
\vspace{0.3cm}
\begin{theorem}\label{thm:convex-hull}
$\!\!\! .\,\,$
For any $\delta, \eps>0$, with probability at least $1 - \delta$,
all of the agents fall within distance~$\eps$ 
of the convex hull of the fixed agents after a number of steps at most 
$$
\frac{1}{p \delta} \left(\frac{c}{\rho} 
                         \log \frac{dn}{\eps}\right)^{2(n-1)}  
$$
for constant $c>0$.
\end{theorem}
\smallskip

\vspace{1cm}
\proof
Let $Q_t$ be the $h$-by-$h$ upper-left submatrix of $P_t$,
where $h= n-k$.
Note that $Q_{\leq t}:= Q_t\cdots Q_0$ coincides with the
$h$-by-$h$ upper-left submatrix $P_{\leq t}$.
Thus, to show that  the mobile agents are attracted to the convex hull
of the fixed ones, it suffices to prove that $Q_{\leq t}$ 
tends to $\mathbf{0}_{h\times h}$.
To do that, we create an averaging system
consisting of $h+1$ agents embedded in $[0,1]$ and
evolving as $y(t+1)= A_t\, y(t)$, where: $y(t)\in \mathbb{R}^{h+1}$;
$y_{h+1}(0)=0$;  $v= (\mathbb{I}_h - Q_t)\mathbf{1}_h$; and  
$$
A_t= \begin{pmatrix}
    Q_t  & v \\
     \mathbf{0}_h^T & 1
\end{pmatrix}.
$$
The system lacks the requisite symmetry to qualify as a general averaging system,
so we again use symmetrization~\cite{chazelle-Energ1-2011} 
by duplicating the $h$ mobile agents
and initializing the embedding of the two copies as mirror-image reflections 
about the origin.  The new evolution matrix is now $\nu$-by-$\nu$, where
$\nu= 2h+1$:
$$
B_t= 
\begin{pmatrix}
    Q_t  & v & \mathbf{0}_{h\times h} \\
                    u & 1-2\|u\|_1  & u \\
     \mathbf{0}_{h\times h}  & v & Q_t 
\end{pmatrix}.
$$
We define the row vector $u\in \mathbb{R}^h$ by setting
its $i$-th coordinate to $\rho$ if $v_i>0$ and $0$ otherwise.
We require that $1- 2\|u\|_1\geq \rho$; hence
$\rho\leq 1/(2d_t+1)$, where $d_t$ is the number of mobile agents
(among the $h$ of them) adjacent in $G_t$ to at least one fixed agent. 
This condition is easily satisfied by setting $\rho\leq 1/2n$.
The evolution follows the update:
$z(t+1)= B_t\, z(t)$, where $z(t)\in [-1,1]^\nu$
and $z_{h+1}(0)=0$.

Let $G^*$ be the augmented $\nu$-vertex graph formed from $G$
and let $G_t^*$ be its subgraph selected at time $t$. Note that,
via $z(t)$, these graphs are embedded in  $[-1,1]^\nu$.
If $\Delta_t$ denotes the length of the longest edge of $G^*$ at time $t$
and $T_\alpha$ is the last time at which the diameter of the system is at least $\alpha$,
then $\Delta_t \geq  \alpha/\nu$ for all $t \leq T_\alpha$ because $G^*$ is connected.
The longest edge in $G^*$ (with ties broken alphabetically) appears in $G_t^*$ 
with probability at least~$p$.
Fix $s\in (0,1]$ and define the random variable $\chi_t$ 
to be $\Delta_t^s$ if the longest edge of $G$ at time $t$ is in 
$G_t$ and $0$ otherwise. By~\cite{chazelle-Energ2-2019}, the maximum $s$-energy 
is at most $\mathcal{E}_s\leq  2^s(3/\rho s)^{\nu-1}$; hence
$$
\mathbb{E}\, T_\alpha
\leq \left(\frac{\nu}{\alpha}\right)^s \mathbb{E} \sum_{t \geq 0} \Delta_t^s
\leq \frac{1}{p} \left(\frac{\nu}{\alpha}\right)^s \mathbb{E} \sum_{t \geq 0} \chi_t
\leq \frac{1}{p} \left(\frac{\nu}{\alpha}\right)^s \mathcal{E}_s
\leq \frac{2}{p} \left(\frac{\nu}{\alpha}\right)^s  \left(\frac{3}{\rho s}\right)^{\nu-1}.
$$
Minimizing the right-hand side over all $s\in (0,1]$ yields
$$
\mathbb{E}\, T_\alpha\leq \frac{6}{p}
                             \left(\frac{3}{\rho} \log \frac{2n-1}{\alpha}\right)^{2(n-1)}.
$$
By Markov's inequality, $\Pr\, \bigl[\, T_\alpha \geq t_\delta \, \bigr] \leq \delta$, where
\begin{equation}\label{t_delta}
t_\delta :=   \frac{6}{p \delta} \left(\frac{3}{\rho} 
                         \log \frac{2n-1}{\alpha}\right)^{2(n-1)}  \, .
\end{equation}
This implies that 
$\|Q_{\leq t} \, \mathbf{1}_h \|_\infty \leq \alpha$, for all $t > t_\delta$,
with probability at least $1-\delta$.
In other words, for any such $t$, it holds that, for $i\leq h$,
$$
q:= \sum_{j=1}^h (P_{\leq t})_{ij}  = \sum_{j=1}^h (Q_{\leq t})_{ij} \leq \alpha.
$$
Trivially, $x_i(t+1)= q\tilde{u} + (1-q)\tilde{v}$, where 
$$
\tilde{u}=  \frac{1}{q}\sum_{j=1}^h (P_{\leq t})_{ij}x_j(0)
\hspace{0.6cm}
\text{and}
\hspace{0.6cm}
\tilde{v}=  \frac{1}{1-q}\sum_{j=h+1}^n (P_{\leq t})_{ij}x_j(0).
$$
Observing that $\tilde{v}$ lies in the convex hull of the fixed agents, we form
the difference $x_i(t+1) - \tilde{v}=  q(\tilde{u}-\tilde{v})$ and note that 
the distance from $x_i(t+1)$ to the hull is bounded
by $q\|\tilde{u}-\tilde{v}\|_2\leq \alpha \sqrt{d}$.
Setting $\alpha = \eps/\sqrt{d}$ completes the proof.
\hfill $\Box$
\proofend

We can extend this result so as to relate convergence to connectivity.
We now produce the random graph $G_t$ from fixed connected $G$
as we did above, but if this results in a graph with more than $m$
connected components, we add random edges picked uniformly from $G$ until
the number of components drops to $m$.
Using Theorem~\ref{s-Energ-Bounds}~(a)   
in the proof above leads to a more refined bound:

\vspace{0.4cm}
\vspace{0.3cm}
\begin{theorem}\label{thm:convex-hull-connect}
$\!\!\! .\,\,$
For any $\delta, \varepsilon > 0$, with probability at least $1 - \delta$, 
all of the agents fall within distance~$\varepsilon$ of the convex hull
of the fixed agents in time bounded by
$$
\frac{c}{p\delta}\left(\frac{1}{\rho}\right)^{2(n-1)}
\left(mn\log \frac{dn}{\varepsilon}\right)^{2m-1},
$$
for constant $c>0$. This assumes that no graph used in the process has more than $m$
connected components.
\end{theorem}
\smallskip
\vspace{0.3cm}

\proof
By Theorem~\ref{s-Energ-Bounds}~(a),  we know that
$\mathcal{E}_{m,s} \leq (b/s)^m(1/\rho)^{n-1}$, for any $s\in (0,1]$, 
where $b=O(mn)$.
The previous proof leads us to update~(\ref{t_delta}) into: 
\begin{equation*}
t_\delta :=   \frac{c}{p \delta} \left(\frac{1}{\rho}\right)^{2(n-1)}
                         \left(mn\log \frac{n}{\alpha}\right)^{2m-1}  \, ,
\end{equation*}
for constant $c>0$, from which the theorem follows.
\hfill $\Box$
\proofend

\newpage


\end{document}